\documentclass[aos,preprint]{imsart}

\usepackage{amsmath}
\usepackage{amssymb}

\def\BibTeX{{\rm B\kern-.05em{\sc i\kern-.025em b}\kern-.08em
    T\kern-.1667em\lower.7ex\hbox{E}\kern-.125emX}}

\numberwithin{equation}{section}
\newtheorem{rmk}{Remark}[section]
\allowdisplaybreaks 
\def\1{\mbox{1\hspace{-.35em}1}} % indicatrice
\def\R{\mathbb{R}}

\def\N{\mathbb{N}}
\def\P{\mathbb{P}}
\def\E{\mathbb{E}}
\def\L{\mathbb{L}}

\def\R{\mathbb{R}}
\def\C{\mathbb{C}}
\def\Z{\mathbb{Z}}

\def\v{\mbox{Var\,}}

\def\lip{\mbox{Lip\,}}
\def\Lip{\mbox{Lip\,}}

\def\limiteproban{\renewcommand{\arraystretch}{0.5}
\begin{array}[t]{c}
\stackrel{{\cal P} }{\longrightarrow} \\
{\scriptstyle n\rightarrow\infty}
\end{array}\renewcommand{\arraystretch}{1}}

\def\limiteloin{\renewcommand{\arraystretch}{0.5}
\begin{array}[t]{c}
\stackrel{{\cal D} }{\longrightarrow} \\
{\scriptstyle n\rightarrow\infty}
\end{array}\renewcommand{\arraystretch}{1}}

\def\limiten{\renewcommand{\arraystretch}{0.5}
\begin{array}[t]{c}
\stackrel{}{\longrightarrow} \\
{\scriptstyle n\rightarrow\infty}
\end{array}\renewcommand{\arraystretch}{1}}

\def\limitepsk{\renewcommand{\arraystretch}{0.5}
\begin{array}[t]{c}
\stackrel{a.s.}{\longrightarrow} \\
{\scriptstyle k\rightarrow\infty}
\end{array}\renewcommand{\arraystretch}{1}}

\def\limiteasn{\renewcommand{\arraystretch}{0.5}
\begin{array}[t]{c}
\stackrel{a.s.}{\longrightarrow} \\
{\scriptstyle n \rightarrow\infty}
\end{array}\renewcommand{\arraystretch}{1}}

\newtheorem{Theorem}{Theorem}

\newtheorem{Proposition}{Proposition}

\newtheorem{Lemma}{Lemma}
\newtheorem{Corollary}{Corollary}
\newtheorem{Remark}{Remark}

\arxiv{math.PR/0000000}

\begin{document}

\begin{frontmatter}

\title{Asymptotic normality of the Quasi Maximum Likelihood Estimator
for multidimensional causal processes}%\protect\thanksref{T1}
\runtitle{Asymptotic normality of the QMLE
for causal processes}
%\thankstext{T1}{Footnote to the title with the `thankstext' command.}

\begin{aug}
  \author{Jean-Marc Bardet \ead[label=e1]{Jean-Marc.Bardet@univ-paris1.fr}
  \ead[label=u1,url]{http://matisse.univ-paris1.fr/bardet/}}
  \and
  \author{Olivier Wintenberger \ead[label=e2]{Olivier.Wintenberger@univ-paris1.fr}
  \ead[label=u2,url]{http://wintenberger.fr}}
  \runauthor{J.-M. Bardet and O. Wintenberger}

  \affiliation{CES (SAMOS-Matisse), University Paris 1 Panth\'eon-Sorbonne, France}

\address{
SAMOS\\
Centre Pierre Mend\`es-France\\
90 Rue de Tolbiac\\
75634 Paris Cedex 13,\\
\printead{e1,e2}\\
          \printead{u1,u2}}

\end{aug}

\begin{abstract}~~Strong consistency and asymptotic normality of the Quasi-Maximum
Likelihood Estimator (QMLE) are given for a general class of
multidimensional causal processes. For particular cases already
studied in the literature (for instance univariate or multivariate
GARCH, ARCH, ARMA-GARCH processes) the assumptions required for
establishing these results are often weaker than existing
conditions. The QMLE asymptotic behavior is also given for
numerous new examples of univariate or multivariate processes (for
instance TARCH or NLARCH processes).
\end{abstract}



\end{frontmatter}

\section{Introduction}
In this paper the asymptotic behavior of the Quasi-Maximum Likelihood Estimator (QMLE) is
studied for general $\R^m$-valued stationary process. The time series
$X=(X_t,~t\in \Z)$ is defined as a solution of the equation:
\begin{equation}\label{eq::sys}
X_t= M_{\theta_0}(X_{t-1},X_{t-2},\ldots) \cdot \xi_t+
f_{\theta_0}(X_{t-1},X_{t-2},\ldots),~~\forall ~t\in\Z,
\end{equation}
almost everywhere (a.e.). Here
$M_{\theta_0}(X_{t-1},X_{t-2},\ldots)$ is a $(m\times p)$-random
matrix having almost surely (a.s.) full rank $m$, the sequence
$(\xi_t)_{t\in\Z}$ of $\R^p$-random vectors $(\xi_t^{(k)})_{1\leq
k \leq p}$ are independent and identically distributed satisfying
$\E\big [\xi_0^{(k)}\xi_0^{(k')}\big ]=0$ for $k\neq k'$ and $\E
\big [{\xi_0^{(k)}}^2\big ]=\v(\xi_0^{(k)})=1$ and
$f_{\theta_0}(X_{t-1},X_{t-2},\ldots)$ is a sequence of $\R^m$-random vectors.
Various popular econometric time series models can be written in the form (\ref{eq::sys}).
The case $f_\theta \equiv 0$ and
\begin{multline}
\label{H}  H_\theta(X_{t-1},X_{t-2},\ldots):=C_0+\sum_{i=1}^{q'}
\sum_{j=1}^k C_{ij}X_{t-i}X'_{t-i}C'_{ij}\\
+\sum_{i=1}^{q} \sum_{j=1}^k D_{ij}H_\theta
(X_{t-i-1},X_{t-i-2},\ldots)D'_{ij}, \mbox{ where\footnotemark } H_\theta
:=M_\theta \cdot M_\theta',
\end{multline}
corresponds to the BEKK representation of multivariate
GARCH($q,q')$ defined by Engle and Kroner \cite{Engle1995},
see also Bollerslev \cite{Bollerslev1990}. Their natural generalization,
\footnotetext{Here $A'$ is the transpose of the matrix $A$.}
$$
H_\theta(X_{t-1},X_{t-2},\ldots):=B_0+\sum_{i=1}^\infty
B_{i}X_{t-i}X'_{t-i}B'_{i},
$$
defines the multivariate ARCH($\infty$) processes.
If $M_\theta \equiv I_d$, a process $X$ satisfying
relation (\ref{eq::sys}) is a multivariate Non Linear AR($\infty$) process.\\

Various methods can be employed to estimate the unknown parameter
$\theta_0$. Maximum Likelihood Estimation (MLE) is a common one.
Several authors studied the asymptotic behavior of MLE for
particular cases of multivariate processes satisfying
(\ref{eq::sys}), see for instance Bollerslev and Wooldridge
\cite{Bollerslev1992}, Jeantheau  \cite{Jeantheau1998} for
multivariate GARCH($q,q'$) processes and Dunsmuir and Hannan
\cite{Dunsmuir1976}, Mauricio \cite{Mauricio1995} for multivariate
ARMA processes. A proof of the efficiency of those estimators was obtained in Berkes and
Horv\'ath \cite{Berkes2004}, in the case of one-dimensional
GARCH($q,q'$). Even if the convergence rate of the MLE can be optimal
this method presents numerous drawbacks. For example, the
conditional likelihood depends on the distribution of the
innovations $\xi_t$, which is often unknown, and on all the
past values of the process $X$, which are unobserved. \\

In the present paper we consider  an approximation of the MLE called
Quasi-Maximum Likelihood Estimation (QMLE). If the sequence
$(\xi_t)_{t\in\Z}$ is a sequence of standardized Gaussian vectors,
the conditional likelihood of $X$ is, up to an additional
constant, equal to
\begin{eqnarray}
\label{eq::MLE}&&L_n(\theta):= -\frac{1}{2}
\sum_{t=1}^nq_t(\theta)\qquad \mbox{for all }\theta\in\Theta\\
 &&\qquad\mbox{with}~~q_t(\theta):=\Big [ \big
(X_t-f^t_\theta\big)' \big ( H_\theta^t \big) ^{-1} \big
(X_t-f^t_\theta\big)+\ \log \big (\det \big (H_\theta^t\big ) \big
) \Big ] \nonumber ,
\end{eqnarray}
$f_\theta^t=f_{\theta}(X_{t-1},X_{t-2},\ldots)$, $M_\theta^t=M_{\theta}(X_{t-1},X_{t-2},\ldots)$ and $H_\theta^t:=M_\theta^t
{M_\theta^t}'.$\\

>From now on we omit any assumption on the distribution of the
$\xi_t$. The QMLE is obtained by plugging in the likelihood the
approximations $ \widehat
f_\theta^t:=f_{\theta}(X_{t-1},\ldots,X_1,u)$, $\widehat
M_\theta^t:=M_{\theta}(X_{t-1},\ldots,X_1,u)$ and
$\displaystyle\widehat H_\theta^t:=\widehat M_\theta^t
\cdot(\widehat M_\theta^t)'$ where $u$ is a finitely-non-zero
sequence\footnote{This means that $u_n\neq0$ only for finitely many $n\in\N$.} $(u_n)_{n\in\N}$:
\begin{eqnarray}\label{eq::QMLE}
~~\widehat L_n(\theta)&\hspace{-3mm}:=&\hspace{-3mm} -\frac{1}{2}
\sum_{t=1}^n \widehat q_t(\theta) \\
\nonumber && \hspace{-1cm} \mbox{with}~~\widehat q_t(\theta):=\Big
[ \big (X_t-\widehat f^t_\theta\big)' \big (\widehat H_\theta^t
\big) ^{-1} \big (X_t-\widehat f^t_\theta\big) + \log \big (\det
\big(\widehat H_\theta^t \big )\big )\Big ].
\end{eqnarray}
The QMLE $\widehat \theta_n$ is the
M-estimator associated with the quasi-likelihood $\widehat L_n$ given as the maximizer
\begin{equation}\label{qmle}
\widehat \theta_n:=\underset{\theta \in \Theta}{\mbox{Argmax }}\widehat
L_n(\theta).
\end{equation}
A basic idea of this paper is to restrict the set of
parameters $\Theta$ in such a way that moment conditions on $\xi_0$ imply both the existence of a solution $X$ and
finite moments of sufficiently high order for $X$. This strategy is available for the very
general model (\ref{eq::sys}) thanks to a result of Doukhan and
Wintenberger \cite{Doukhand}, see Section \ref{Existence}. Then we
use the moment conditions to settle both consistency and
asymptotic normality, see Section \ref{TLC}.\\

We restrict the set of the parameters in such a way that we only assume
finite moments of orders $2$ or $4$ on $\xi_0$, which are
necessary conditions for consistency or asymptotic normality,
respectively, see for example Straumann and Mikosch \cite{Straumann2006} for some particular classes of non-linear time series models. In turn, these conditions guarantee
the existence of moments of order $2$ or $4$ of $X$, respectively. Notice that for one-dimensional GARCH models  these moment conditions on $X$ can be relaxed, see Francq and Zako\"{\i}an \cite{Francq2004},
Berkes {\it et al.} \cite{Berkes2003}. For Markovian models,
Straumann and Mikosch \cite{Straumann2006} achieved the asymptotic
normality assuming moment conditions but the corresponding
restriction on $\Theta$ is non-explicit except for the AGARCH
models. In the case of ARCH($\infty$), the conditions are not
comparable with those in Robinson and Zaffaroni
\cite{Robinson2006}. Our restriction on $\Theta$ is  stronger
whereas we sharpen the moment conditions of order $2+\delta$
to the order $2$ on $\xi_0$ for the strong consistency. Finally,
for multivariate models the conditions are sharper than those in Comte
and Lieberman \cite{Comte2003} and Ling and McAleer \cite{Ling2003}
who derived the asymptotic normality for particular models under moments of order $4$,
$6$ or $8$ on $X$. In Section \ref{Examples} we provide for the first time the consistency and asymptotic normality
of the QMLE in TARCH, NLARCH and Non Linear AR($\infty$) models.\\

But to begin with, the following Section \ref{Existence} deals
with the various assumptions on the general model (\ref{eq::sys})
that are needed.
\section{Notation and assumptions} \label{Existence}
In the sequel, some standard notation is used:
\begin{itemize}
\item The symbol $\| . \|$ denotes the usual Euclidean norm
of a vector or a matrix (for $A$ a $(n \times p)$-matrix,
$\|A\|=\sup_{\|Y\|\leq 1} \big \{\|AY\|,~Y\in \R^p \big \}$);
\item For the measurable
vector- or matrix-valued function $g$ defined on $\Theta$, $\| g
\|_\Theta=\sup_{\theta \in \Theta} \| g(\theta) \|$;
\item If $V$ is a vector space then $V^{\infty}$ denotes the set of the finitely-non-zero sequences $x$ i.e.,  there exists $N>0$ such that $x=(x_1,x_2,\ldots,x_N,0,0,\ldots)$;
\item  The symbol $0$ denotes the null sequence in $\R^\N$;
\item If $V$ is a Banach space and $\Theta$ is a subset of $\R^d$ then ${\cal C}(\Theta,V)$ denotes the
Banach space of $V$-valued continuous functions on $\Theta$
equipped with the uniform norm $\| \cdot \|_\Theta$ and $\L^r({\cal
C}(\Theta,V))$ ($r\ge 1$) denotes the Banach space of random a.e. continuous
functions $f$ such that $\E\big [\|f\|_\theta^r\big ]<\infty$.
\end{itemize}
%
% When it exists, the limits of $f_\theta(x_1,\ldots,x_t,0,\ldots)$ and $g_\theta(x_1,\ldots,x_t,0,\ldots)$ as $t\to \infty$ are denoted $f_\theta(x_1,x_2,\ldots,x_t,\ldots)$ and $g_\theta(x_1,x_2,\ldots,x_t,\ldots)$. Assume moreover that:\\
\subsection{Definition of the parameter sets $\Theta(r)$ and $\widetilde \Theta(r)$}
In proposition \ref{stationarity} below we provide the existence of a stationary solution of the general model \eqref{eq::sys}. Two conditions of different types are used: the first one is a Lipschitz condition on the functions $f$ and $M$ in \eqref{eq::sys}, the second one is a restriction on the set of the parameters.\\

Let us assume that for any $\theta\in\R^d$, $x\mapsto f_\theta(x)$ and $x\mapsto M_\theta (x)$ are Borel functions on $({\R^m}) ^\infty$ and that $\mbox{Rank }M_\theta (x)= m$ for all $x\in({\R^m}) ^\infty$. Assume that there exist two sequences $(\alpha_j(f,\theta))_{j\ge1}$ and $(\alpha_j(M,\theta))_{j\ge1}$ satisfying, for all $x$,
$y$ in $({\R^m}) ^{\infty}$,
$$\left \{
\begin{array}{lll}
\|f_{\theta}(x)-f_{\theta}(y)\|
&\le &
\sum_{j=1}^\infty \alpha_j(f,\theta)\|x_j-y_j\|,\\
\|M_{\theta}(x)-M_{\theta}(y)\|
&\le &\sum_{j=1}^\infty \alpha_j(M,\theta)\|x_j-y_j\|.
\end{array}
\right .
$$
For some models, as mentioned in remark \ref{tildeh}, it can be more efficient to replace the condition on $M$ by the existence of a
sequence $(\alpha_j(H,\theta))_{j\ge1}$ such that
$$
\|H_{\theta}(x)-H_{\theta}(y)\|\le \sum_{j=1}^\infty \alpha_j(H,\theta)\|x_jx'_j-y_jy'_j\|,
$$
where $H_\theta
:=M_\theta \cdot M_\theta'$. Assuming $\E \|\xi_0\|^r<+\infty$ for some $r>0$, we can define the set
\begin{equation}
\label{propINF}\Theta(r)=\left\{\theta\in\R^d~\Big/~ \sum_{j=1}^\infty
\alpha_j(f,\theta)+\left(\E\|\xi_0\|^r\right)^{1/r}\sum_{j=1}^{\infty} \alpha_j(M,\theta)<1\right\}.
\end{equation}
This set depends on the distribution of $\xi_0$ via the moments $\E\|\xi_0\|^r$. But thanks to the fact that $\E\big [\xi_0^{(k)}\xi_0^{(k')}\big ]=0$ for $k\neq k'$ and $\E
\big [{\xi_0^{(k)}}^2\big ]=\v(\xi_0^{(k)})=1$ the set $\Theta(2)$ simplifies:
$$
\Theta(2)=\left\{\theta\in\R^d~\Big/~  \sum_{j=1}^\infty
\alpha_j(f,\theta)+\sqrt{p}\sum_{j=1}^{\infty} \alpha_j(M,\theta)<1\right\}.
$$
\begin{Proposition}\label{stationarity}
If $\theta_0\in \Theta(r)$ for some $r\ge1$ there exists a unique
causal ($X_t$ is independent  of $(\xi_i)_{i> t}$ for $t\in\Z$) solution $X$ to the equation \eqref{eq::sys} which is
stationary and ergodic  and satisfies
$\E\big \|X_0\big\|^r <\infty$.
\end{Proposition}
This result generalizes the one proved by Giraitis {\it et al.}
\cite{Giraitis2000} for ARCH($\infty$) models. It automatically yields weak dependence properties, see \cite{Doukhand} for details. For such non Markovian models, the classical Lyapunov condition of Bougerol \cite{Bougerol1993} cannot be applied.\\

Let us now consider the special cases of \eqref{eq::sys} where $f\equiv 0$, $m=p=1$ and there exists a Borel function $\widetilde H_\theta$ such that $H_\theta(x)=\widetilde H_\theta(x^2)$ for all $x\in \R^\infty$.
\begin{Corollary}\label{scal}The result of Proposition \ref{stationarity} holds if $\theta_0\in\widetilde \Theta(r)$ for $r\ge2$ where
\begin{equation}
\label{CondH}\widetilde \Theta(r)=\left\{\theta\in\R^d~\Big/~ \E |\xi_0|^r
\big(\sum_{j=1}^\infty \alpha_j(H,\theta)\big)^{r/2}<1\right\}.
\end{equation}
\end{Corollary}
\begin{rmk}\label{tildeh}
\normalfont The ARCH($\infty$) process was defined by Robinson
\cite{Robinson1991} as solution of the model:
\begin{eqnarray}\label{arch}
X_t = \sigma_t \xi_t, \qquad \sigma^2_t = b_0(\theta_0) +
\underset{j = 1}{\overset{\infty}{\sum}} b_j(\theta_0) X^2_{t - j},
\end{eqnarray}
where, for all $\theta\in\R^d$, $(b_j(\theta))_{j\ge1 }$ are sequences
of non-negative real numbers. Here, $f\equiv 0$, $p=m=1$, $\alpha_j(M,\theta)=\sqrt{b_j(\theta)}$ and $\alpha_j(H,\theta)=b_j(\theta)$. Working with the set $\widetilde \Theta(r)$, larger than $\Theta(r)$, gives more general results.
\end{rmk}
\subsection{Uniform assumptions on $\Theta$}
Fix some compact subset $\Theta$ of $\R^d$. For  any
sequences $x$, $y$ of $({\R^m}) ^{\infty}$, the functions $\theta
\mapsto f_{\theta}(x)$ and $\theta \mapsto M_\theta(x)$ are
assumed to be continuous on $\Theta$. As in \cite{Straumann2006},
uniform continuity conditions on $\Theta$ are required to apply
the QMLE procedure, see Lemma \ref{lem::estvar} of the Section
\ref{TLC}. Assume that $\|f_\theta(0)\|_{\Theta}<\infty$ and $
\|M_\theta(0)\|_{\Theta} <\infty$. To settle the
assumptions in a short way, let us introduce the generic symbol $\Psi$ for any of the functions
$f$, $M$ or $H$.
\begin{description}
\item[(A1($\Psi$))] Let $\alpha_j(\Psi)=\sup_{\theta\in\Theta}\alpha_j(\Psi,\theta)$
be such that $\sum_{j\ge1} \alpha_j(\Psi)<\infty$.
\item[(A2)] There exists $\underline H>0$ such that
$\inf_{\theta\in\Theta}\det \big(H_\theta(x) \big ) \ge
\underline H$ for all $x\in ({\R^m}) ^{\infty}$.
\item[(A3($\Psi$))] The function $\theta \in \Theta  \mapsto
\Psi_\theta(x)$ is $2$ times continuously differentiable for all
$x\in(\R^m)^{\infty}$ and
$$ \Big \|\frac {\partial
\Psi_\theta(0)}{\partial \theta} \Big \|_\Theta+\Big \|\frac
{\partial^2 \Psi_\theta(0)}{\partial \theta\partial \theta'} \Big
\|_\Theta< \infty.$$
Moreover assume that there exist two integrable sequences
$\big(\alpha^{(i)}_j(\Psi)\big)_{j\ge1}$, $i=1,2$, such that for all $x$, $y\in(\R^m)^{\infty}$
\begin{eqnarray*}\Big \|\frac
{\partial \Psi_\theta(x)}{\partial \theta
 } -\frac {\partial \Psi_\theta(y)}{\partial \theta} \Big \|_\Theta &\le& \displaystyle  \sum_{j=1}^\infty
\alpha^{(1)}_j(\Psi)
\|x_j-y_j\|,\\
\Big \|\frac
{\partial^2 \Psi_\theta(x)}{\partial \theta\partial \theta'
 } -\frac {\partial^2 \Psi_\theta(y)}{\partial \theta\partial \theta'} \Big \|_\Theta &\le& \displaystyle  \sum_{j=1}^\infty
\alpha^{(2)}_j(\Psi)
\|x_j-y_j\|.
\end{eqnarray*}
If $\Psi=H$, $\|x_j-y_j\|$ in the RHS terms is replaced with $\|x_jx_j'-y_jy_j'\|$.
\end{description}
The last assumption on the derivatives is just needed for the asymptotic normality of the QMLE.

\subsection{Identifiability and variance conditions}
We assume the same identifiability condition as in Jeantheau \cite{Jeantheau1998}:
\begin{description}
\item[(Id)] For all $\theta\in \Theta$,
($f^t_\theta=f^t_{\theta_0}$ and
$H^t_\theta=H^t_{\theta_0}$ a.s.) $\Rightarrow\theta=\theta_0$.\\
\item[(Var)] One of the families $({\partial f_{\theta_0}^t}/{\partial
\theta_i})_{1\le i\le d}$ or $({\partial
H_{\theta_0}^t}/{\partial \theta_i})_{1\le i\le d}$ is a.e.
linearly independent, where:
$$
\frac{\partial f_\theta^t}{\partial \theta}:=
\frac{\partial f_\theta}{\partial \theta}(X_{t-1},\ldots)\mbox{ and }
\frac{\partial H_\theta^t}{\partial \theta}:=
\frac{\partial H_\theta}{\partial \theta}(X_{t-1},\ldots).
$$
\end{description}
The condition {\bf (Var)} is needed for ensuring finiteness of the
asymptotic variance in the result on asymptotic normality. For ARCH($\infty$),
Robinson and Zaffaroni \cite{Robinson2006} give sufficient
assumptions for both {\bf (Id)} and {\bf (Var)}. They are easier to verify than {\bf (Id)} and {\bf
(Var)} but are not as general. Alternative conditions similar to those for ARCH($\infty$) are
not straightforward in the general model \eqref{eq::sys} because of its
non-linear character.

\section{Asymptotic behavior of the QMLE} \label{TLC}
If the model satisfies the conditions of Corollary \ref{scal}, the
set $\Theta(r)$ can be replaced with $\widetilde\Theta(r)$ in all
the results of this section.
\subsection{Invertibility}
Here we follow the presentation of Straumann and Mikosch
\cite{Straumann2006}. The approach of the QMLE is based on an
approximation of $f^t_\theta=\E(X_t~|~ X_{t-1},X_{t-2},\ldots)$
and $H^t_\theta =\E\big((X_t-f^t_\theta)(X_t-f^t_\theta)'~|~
X_{t-1},X_{t-2},\ldots\big)$ by $\widehat f^t_\theta$ and
$\widehat H^t_\theta $, defined as in the introduction. Invertibility is the property that $\widehat f^t_\theta$ and $\widehat
H^t_\theta $ converge to the unobservable $f^t_\theta$ and
$H^t_\theta$, see Section 3.2 of \cite{Straumann2006} for more
details. The following lemma states this result which is a
necessary step in the proof of the QMLE consistency.
\begin{Lemma}\label{lem::estvar}
Assume that $\theta_0\in \Theta(r)$ for $r \geq 2$ and that $X$ is the stationary
solution of the equation \eqref{eq::sys}.
\begin{enumerate}
\item If {\bf (A1(f))} holds then $f_\theta^t\in \L^r({\cal C}(\Theta,\R^m))$ and
\begin{equation}\label{eq::estf}
\E \big [\|\widehat f_\theta^t-f_\theta^t\|_{\Theta}^r \big ] \le
\E\big [\|X_0\|^r \big ]\Big (\sum_{j\ge t}\alpha_j(f)\Big )^r
~~\mbox{for all $t \in \N^*$}.
\end{equation}
\item If {\bf (A1(M))} holds then $H_\theta^t\in \L^{r/2}({\cal C}(\Theta,{\cal M}_m))$ and there exists $C>0$ not depending on $t$ such that
\begin{equation}\label{eq::estfM}
\E \big [\|\widehat H_\theta^t-H_\theta^t\|^{r/2}_\Theta \big]\le
C \Big ( \sum_{j\ge t}\alpha_j(M)\Big )^{r/2} ~~\mbox{for all $t
\in \N^*$}.
\end{equation}
\item If {\bf (A1(H))} holds then $H_\theta^t\in \L^{r/2}({\cal C}(\Theta,{\cal M}_m))$ and
\begin{equation}\label{eq::estH}
\E \big [\|\widehat H_\theta^t-H_\theta^t\|_{\Theta}^{r/2} \big ] \le
\E\big [\|X_0\|^{r} \big ]\Big (\sum_{j\ge t}\alpha_j(H)\Big )^{r/2}
~~\mbox{for all $t \in \N^*$}.
\end{equation}
\end{enumerate}
Moreover, under any of the two last conditions and with {\bf (A2)}, $ H_\theta^t$ is an
invertible matrix and $\Big\|\big (\widehat H_\theta^t  \big)
^{-1}\Big\|_\Theta\le \underline H^{-1/m}$.
\end{Lemma}
The proof is given in Section \ref{proof1}.
\subsection{Strong consistency}
In the following theorem, we assume by convention that if {\bf (A1(M))} holds then
$\alpha_j(H)=0$ and if {\bf (A1(H))} holds then $\alpha_j(M)=0$.
\begin{Theorem}\label{th::consist}
Assume that $\theta_0\in \Theta$ for a compact subset $\Theta\subset\Theta(2)$. Let $X$ be the stationary
solution of the equation \eqref{eq::sys}. Let {\bf (A1(f))}, {\bf (A2)} and {\bf (Id)} hold.
Moreover, if {\bf (A1(M))} or {\bf (A1(f))} hold with
\begin{eqnarray}\label{condP}
\alpha_j(f)+\alpha_j(M)+\alpha_j(H)=O\big (j^{-\ell}\big )~~\mbox{for some}~~\ell
>3/2,
\end{eqnarray}
then the QMLE $\widehat \theta_n$ defined by
\eqref{qmle} is strongly consistent, {\it i.e.}
$
\widehat \theta_n \limiteasn \theta_0.
$
\end{Theorem}
The proof is given in Section \ref{proof2}.
\subsection{Asymptotic normality}
We use the following convention: if {\bf (A3(M))} holds then
$\alpha_j^{(1)}(H)=0$ and if {\bf (A3(H))} holds then
$\alpha_j^{(1)}(M)=0$.
\begin{Theorem}\label{th::asnorm}
Assume that $\theta\in \stackrel{\circ}{ \Theta}$ with
$\stackrel{\circ}{ \Theta}\subset \Theta(4)$ where
$\stackrel{\circ}{ \Theta}$ denotes the interior of a compact
subset $\Theta\subset\R^d$. Let $X$ be the stationary solution of
the equation \eqref{eq::sys}. Assume that the conditions of
Theorem \ref{th::consist} and {\bf (A3(f))}, {\bf (Var)}
hold. Moreover, if {\bf (A3(M))} or {\bf (A3(H))} holds with
\begin{equation}\label{condth2}
\alpha^{(1)}_j(f)+\alpha^{(1)}_j(M)+\alpha^{(1)}_j(H)=O\big (j^{-\ell'}\big
)~~\mbox{for some}~~\ell' >3/2,
\end{equation}
then the QMLE $\widehat\theta_n$ is strongly consistent and asymptotically
normal, i.e.,
\begin{eqnarray}\label{tlcqmle}
\sqrt{n}\big (\widehat\theta_n-\theta_0\big )\limiteloin {\cal
N}_d\big (0 \ , \ F(\theta_0)^{-1} G(\theta_0)
F(\theta_0)^{-1}\big ),
\end{eqnarray}
where the matrices $F(\theta_0)$ and $G(\theta_0)$ are
defined in \eqref{F0} and \eqref{G0} respectively.
\end{Theorem}
The proof is given in section \ref{proof3}.

\section{Examples}\label{Examples}
In this section, the previous asymptotic results are applied to
several examples. For ARCH, GARCH, AR and GARCH-ARMA processes,
the consistency and asymptotic normality have already been
settled and we compare the different conditions from the literature with ours. For
other examples, such as TARCH, multivariate ARCH and NLARCH processes,
the consistency and the asymptotic normality of the QMLE are novel
results. Examples satisfying the conditions of Corollary
\ref{scal} are studied first.
%Traditional GARCH($q,q'$) or ARCH($\infty$) processes where
%$m=p=1$ are particular cases of processes satisfying
%(\ref{eq::sys}).

\subsection{{\em ARCH($\infty$)} processes}
By Remark \ref{tildeh}, the set $\widetilde \Theta(r)$ is
well-adapted to that case
\begin{equation}\label{archstbis}
\widetilde \Theta(r)=\left\{\theta\in\R^d~\Big/~ \sum_{j=1}^\infty b_j(\theta) <\Big ( \E
\big [ | \xi_0|^r \big ] \Big )^{-2/r}\right\}.
\end{equation}
For $\theta_0\in \widetilde\Theta(r)$, the existence of a
stationary solution and of its $r$-th order moments is also
settled in Giraitis {\it et al.} \cite{Giraitis2000}. For an
excellent survey about results and applications of ARCH models, we
refer the reader to Giraitis {\it et al.} \cite{Giraitis2006}. Here we formulate a version of Theorems \ref{th::consist} and
\ref{th::asnorm} adapted to the context.
\begin{Proposition}\label{QMLEarchuniv}
Let $\Theta$ be a compact subset of $\widetilde\Theta(2)$ and $X$ the stationary solution of
\eqref{arch}. Assume that $\inf_{\theta \in \Theta} b_0(\theta)>0$ and
that $\theta\mapsto b_j(\theta)$ be continuous functions satisfying
$$
\sup_{\theta \in \Theta}~b_j(\theta) =O\big (j^{-\ell}\big
)~~\mbox{for some}~~\ell >3/2.
$$
\begin{enumerate}
\item If {\bf (Id)} holds then the QMLE
$\widehat\theta_n$ is strongly consistent.
\item Assume that $\theta_0\in  \stackrel{\circ}{ \Theta}$ with
$\stackrel{\circ}{ \Theta}\subset \widetilde\Theta(4)$, that
$\xi_0^2$ has a non-degenerate distribution. Let the functions $\theta
\mapsto b_j(\theta)$ be $2$-times continuously
differentiable on $\Theta$ for all $j \in \N$ satisfying for all $(k,k')\in
\{1,\ldots,d\}^2$,
$$
\sup_{\theta \in \Theta}\Big | \frac {\partial b_j(\theta)}
{\partial \theta_k } \Big | =O\big (j^{-\ell'}\big
)~~\mbox{for some}~~\ell'>3/2\\
\mbox{ and }\sum_{g\ge1}\sup_{\theta \in \Theta}
\Big | \frac {\partial^2 b_j(\theta)} {\partial \theta_k
\partial \theta_{k'}} \Big |<\infty;
$$
If there exists a finite subset $A \subset
\N$ such that
\begin{equation}\label{asrob}
 \Big ( \Big (\frac{\partial
b_j(\theta_0)}{\partial \theta_k} \Big )_{j \in A} \Big )_{1\le k\le
d}\mbox{ is linearly independent,}
\end{equation}
then the QMLE
$\widehat\theta_n$ is asymptotically normal, {\it i.e.}, it satisfies
\eqref{tlcqmle}.
\end{enumerate}
\end{Proposition}
For the asymptotic normality of the QMLE we use the condition in equation
(\ref{asrob}) coming from Robinson and Zaffaroni
\cite{Robinson2006} that ensures both {\bf (Id)} and {\bf (Var)}.
Let us compare the results of Proposition \ref{QMLEarchuniv} with those of Theorems 1 and 2
in Robinson and Zaffaroni \cite{Robinson2006}. Those authors obtained the almost sure convergence of the QMLE under
moments of order $r>2$ (instead of $r=2$ here) and a decreasing
rate $j^{-\ell}$ with $\ell>1$ (instead of $\ell>3/2$ here) for
the sequence $(\sup_{\theta\in\Theta}|b_j(\theta)|)_{j\ge 1}$. Concerning the asymptotic
normality for $r=4$, their conditions on both the first
derivatives of $\theta \mapsto b_j(\theta)$ are the same as in
Proposition \ref{QMLEarchuniv}. They required also conditions on the third
derivatives (nothing like this here).
%For the specific parametric choice $\theta=(A,B,\alpha,\beta)$ and
%$b_j(\theta)=A\alpha^j+Bj^{-\beta}$ and $\Theta$ in order to
%derive ST$(r)$ with $A\ge \underline H$, $\beta>1$, $0<\alpha<1$
%and $B\ge0$, the assumptions D$^{(1)}$ and D$^{(2)}$ are
%satisfied.

\subsection{{\em GARCH($q,q'$)} models}
The GARCH($q,q'$) models have been introduced by Engle
\cite{Engle1982}. Here $X$ is the stationary solution
of
\begin{eqnarray}\label{garch}
X_t = \sigma_t \xi_t, \qquad \sigma^2_t =c_0(\theta_0) +
\underset{j = 1}{\overset{q}{\sum}} c_j(\theta_0) X^2_{t - j}+
\underset{j = 1}{\overset{q'}{\sum}} d_j(\theta_0) \sigma^2_{t -
j},
\end{eqnarray}
where $c_j(\theta)$ and $d_j(\theta)$ are non negative real
numbers for all $\theta\in\Theta$. This model can be embedded in the
class of ARCH($\infty$) models (see Giraitis {\it et al.}
\cite{Giraitis2006}), as one needs to set for all $z\in \C$
\begin{equation}\label{condgarch}
b_0(\theta) := \frac{c_0(\theta)}{1 - \sum_{j=1}^{q'}
d_j(\theta)}~~\mbox{and}~~\sum_{i=1}^\infty b_i(\theta) z^i :=
\frac{\sum_{i=1}^q c_i(\theta) z^i}{1 - \sum_{i=1}^{q'}
d_i(\theta) z^i}.
\end{equation}
In the last formula, both the polynomials are supposed to be
coprime.  The results of Theorems \ref{th::consist} and
\ref{th::asnorm} lead to the consistency and asymptotic normality of the QMLE in that case. However our conditions are not as sharp as those in Berkes and Horvath \cite{Berkes2004},
Francq and Zako\"{\i}an \cite{Francq2004} or Straumann and
Mikosch \cite{Straumann2006}.
As a consequence of the expression \eqref{condgarch}, the sequence $(\sup_{\theta \in \Theta}
b_j(\theta))_j$ decreases exponentially fast to $0$ and {\bf
A1(H)} holds automatically. From Corollary \ref{scal}, if $\theta_0\in\widetilde
\Theta(r)$, where $\widetilde \Theta(r)$ is defined as in
(\ref{archstbis}), the GARCH process has solutions of order $r$ as
for ARCH($\infty$). For instance, if $q=q'=1$ and
$\theta=(c_0,c_1,d_1)$, we achieve the optimal condition of
existence of a second-order stationary solution as
$\widetilde\Theta(2)=\{c_1,d_1$ such that $c_1>0$ and $c_1+d_1<1\}$. In the property below, we use
the identification condition of Francq and Zako\"{\i}an~\cite{Francq2004}.
\begin{Proposition}\label{QMLEgarch}
Assume that $\Theta$ is a compact subset of $\widetilde\Theta(2)$
defined in \eqref{archstbis} and that $X$ is the stationary solution
to \eqref{garch}. Assume that $\inf_{\theta\in \Theta}
c_0(\theta)>0$, that $\theta\mapsto c_j(\theta)$ and $\theta\mapsto
d_j(\theta)$ are continuous functions injective on $\Theta$ for all
$j$. If $\xi_0$ has a non degenerate distribution then the QMLE
$\widehat\theta_n$ is strongly consistent.\\
Moreover if $\theta_0\in  \stackrel{\circ}{ \Theta}$ with
$\stackrel{\circ}{ \Theta}\subset \widetilde\Theta(4)$, the
functions $\theta \mapsto c_j(\theta)$  and $\Theta  \mapsto
d_j(\theta)$ are $2$ times continuously differentiable on $\Theta$
satisfying for all $(k,k')\in \{1,\ldots,d\}^2$,
$$
\sup_{\theta \in \Theta}\Big | \frac {\partial b_j(\theta)}
{\partial \theta_k } \Big | =O\big (j^{-\ell'}\big
)~~\mbox{for some}~~\ell'>3/2 \mbox{ and }\sum_{j\ge1}\sup_{\theta \in
\Theta} \Big | \frac {\partial^2 b_j(\theta)} {\partial \theta_k
\partial \theta_{k'}} \Big |<\infty.
$$
then the QMLE $\widehat\theta_n$ is also asymptotically normal.
\end{Proposition}
Our assumptions induce the finiteness of second order moments of $X$.
However Jeantheau \cite{Jeantheau1998} proves that moment conditions for models satisfying the Markov property are not
needed for the consistency of the QMLE. In the case
$\theta=(c_0,c_1,\ldots,c_q,d_1,\ldots,d_{q'})$ the Proposition
\ref{QMLEgarch} simplifies:
\begin{Proposition}\label{QMLEgarchbis}
Assume that $\Theta$ is a compact subset of $\widetilde\Theta(2)$ and that $X$ is the solution of
\eqref{garch}. Then, if $\xi_0^2$ has a non-degenerate distribution,
the QMLE
$\widehat\theta_n$ is strongly consistent.
Moreover if $\theta_0\in  \stackrel{\circ}{ \Theta}$ with
$\stackrel{\circ}{ \Theta}\subset \widetilde\Theta(4)$, then
$\widehat\theta_n$ is also asymptotically normal.
\end{Proposition}

\subsection{{\em TARCH($\infty$)} models}
The process $X$ is called
Threshold ARCH($\infty$) if it satisfies the equations
\begin{multline}
\label{tarch} X_t = \sigma_t \xi_t, \\
\sigma_t = b_0(\theta_0) + \underset{j = 1}{\overset{\infty}{\sum}}
\Big [b_j^+(\theta_0) \max ( X_{t - j},0)- b_j^-(\theta_0) \min
(X_{t - j},0) \Big ],
\end{multline}
where the parameters $b_0(\theta)$, $b_j^+(\theta)$ and
$b_j^-(\theta)$ are assumed to be non negative real numbers. This
class of processes is a generalization of the class of
TGARCH($p$,$q$) processes (introduced by Rabemananjara and Zako\"{\i}an
\cite{Rabemananjara1993}) and AGARCH($p$,$q$) processes
(introduced by Ding {\it et al.} \cite{Ding1993}). Here
$$
\Theta(r)=\left\{\theta\in\R^d~\Big/~ \sum_{j=1}^\infty \max \big
(b_j^-(\theta),b_j^+(\theta)\big ) \leq \Big ( \E \big [ |
\xi_0|^r \big ] \Big )^{-1/r}\right\}
$$
since $\alpha_j(M,\theta)= \max \big
(b_j^-(\theta),b_j^+(\theta)\big )$. Consequently, we can settle for the first time the strong consistency and asymptotic normality of the QMLE for TARCH($\infty$) models:
\begin{Proposition}\label{QMLEtarch}
Let $\Theta$ be a compact subset of $\Theta(2)$, $X$ be the stationary solution to
\eqref{tarch} and assume that {\bf (Id)} holds. Assume that $\inf_{\theta\in \Theta} b_0(\theta)>0$  and
$$
\sup_{\theta \in \Theta}\max \big (b_j^-(\theta),b_j^+(\theta)\big
)  =O\big (j^{-\ell}\big )~~\mbox{for some}~~\ell >3/2,
$$
then the QMLE is strongly consistent.\\
Moreover if $\theta_0\in  \stackrel{\circ}{ \Theta}$ with
$\stackrel{\circ}{ \Theta}\subset \widetilde\Theta(4)$, assume that the
functions $\theta\mapsto b_0(\theta)$, $\theta\mapsto b^+_j(\theta)$
and $\theta \mapsto b^-_j(\theta)$ are $2$ times continuously
differentiable on $\Theta$  satisfying,
\begin{multline*}
\sup_{\theta \in \Theta}\max \Big ( \Big | \frac {\partial
b^+_j(\theta)} {\partial \theta_k } \Big |\, , \, \Big | \frac
{\partial b^+_-(\theta)} {\partial \theta_k } \Big | \Big ) =O\big
(j^{-\ell'}\big
)~~\mbox{for some}~~\ell'>3/2~~
\mbox{and}\\
\sum_{j\ge1}\sup_{\theta \in \Theta}\max\Big ( \Big | \frac
{\partial^2 b^+_j(\theta)} {\partial \theta_k \partial \theta_{k'}}
\Big |\, , \, \Big | \frac {\partial^2 b^-_j(\theta)} {\partial
\theta_k
\partial \theta_{k'}} \Big | \Big )<\infty\mbox{  for all }(k,k')\in
\{1,\ldots,d\}^2.
\end{multline*}
If {\bf (Var)} holds then the QMLE
$\widehat\theta_n$ is also asymptotically normal.
\end{Proposition}

\subsection{Multivariate {\em ARCH($\infty$)} processes}
The multivariate ARCH($\infty$) processes are defined as solutions to equation
(\ref{eq::sys}) where
%with a symmetric
%$M_\theta=H_\theta^{1/2}$ (and therefore $M_\theta
%(M_\theta)'=H_\theta$) and
\begin{eqnarray}\label{BEKK}
H_\theta(X_{t-1},X_{t-2},\ldots):=B_0(\theta)+\sum_{i=1}^\infty
B_{i}(\theta)X_{t-i}X'_{t-i}B'_{i}(\theta).
\end{eqnarray}
Here $B_{i}(\theta)$ is a non-negative definite
$d \times d$ matrice for all $\theta\in\R^d$. As $\alpha_j(M,\theta)=B_j(\theta)$
we have
\begin{eqnarray}\label{condARCH}
\Theta(r)=\left\{\theta\in\R^d~\Big/~\sum_{j=1}^\infty B_j(\theta)<  \Big ( \E \big [ \| \xi_0\|^r \big
] \Big )^{-1/r}\right\}.
\end{eqnarray}
\begin{Proposition}\label{QMLEarch}
Let $\Theta$ be a compact subset of $\Theta(2)$ and $X$ be the stationary solution to \eqref{eq::sys} when relation \eqref{BEKK} holds.
Assume that $\inf_{\theta \in \Theta} \det B_0(\theta)
>0$, {\bf (Id)} holds and
$$
\| B_j(\theta)\|_\theta =O\big
(j^{-\ell}\big )~~\mbox{for some}~~\ell >3/2.
$$
Then the QMLE is strongly consistent.\\
Moreover, if $\theta_0\in  \stackrel{\circ}{ \Theta}$ with
$\stackrel{\circ}{ \Theta}\subset \widetilde\Theta(4)$, assume that
the functions $\theta \mapsto B_j(\theta)$ are $2$ times
continuously differentiable on $\Theta$ satisfying for all
$(k,k')\in \{1,\ldots,d\}^2$,
$$
\Big \| \frac {\partial B_j(\theta)} {\partial \theta_k } \Big
\|_\Theta= O\big (j^{-\ell'}\big )~~\mbox{for some}~~\ell'>3/2~~
\mbox{and }\sum_{j\ge1}\Big \| \frac {\partial^2 B_j(\theta)}
{\partial \theta_k \partial \theta_{k'}} \Big \|_\Theta<\infty.
$$
If {\bf (Var)} holds, then the QMLE $\widehat\theta_n$ is also
asymptotically normal.
\end{Proposition}
For the best of our knowledge, the asymptotic behavior of the QMLE for
such models is studied here for the first time.
\subsection{Multivariate {\em GARCH($q,q'$)} models}
Multivariate GARCH($q,q'$) models refer classically to both VEC
and  BEKK models. We refer the reader to Section \ref{armagarch} for VEC models which are subcases of ARMA-GARCH models. BEKK
processes are solutions of equation (\ref{BEKK}) or equivalently
$$
vec(H_t)=vec(C_0)+\sum_{i=1}^qC_i(\theta_0)^\ast
vec(X_{t-i}X_{t-i}')+\sum_{i=1}^pD_i(\theta_0)^\ast vec(H_{t-i}),
$$
where $vec$ is the operator that stacks together the columns of a matrix. For any $p\times k$ matrix $A$:
$A_i^\ast=\sum_{j=1}^k A_{i,j}\otimes A_{i,j}$ for $i=1,\ldots,p$, where $\otimes$
denoting the Kronecker product.
The multivariate ARCH($\infty$) representation holds with $B_j$ satisfying
\begin{multline}\label{garchbekk}B_0^\ast:= \Big(1 - \sum_{j=1}^{q'}
D_j^\ast\Big)^{-1}\hspace{-5mm} \times C_0^\ast~~\mbox{
and}\\\sum_{i=1}^\infty B_i^\ast Z^i := \Big(1 - \sum_{i=1}^{q'}
D_i^\ast Z^i\Big)^{-1}\hspace{-5mm} \times\sum_{i=1}^q C_i^\ast
Z^i~~\mbox{for all}~~Z\in \C^m.
\end{multline}
In the last formula, both the polynomials are supposed to be
coprime.\\

The natural choice $\theta=(C_0,C_1,\ldots,C_q,D_1,\ldots,D_{q'})$
implies that {\bf (A1(M))} and {\bf (Var)} are satisfied. Using
the identification condition of Comte and Lieberman
\cite{Comte2003}, Proposition
\ref{QMLEarch} becomes more simple:
\begin{Proposition}\label{QMLEgarchmult}
Let $\Theta$ be a compact subset of $\Theta(2)$ defined in
\eqref{condARCH} and $X$ be the stationary solution to \eqref{eq::sys} when relations
\eqref{BEKK} and \eqref{garchbekk} hold. If $\inf_{\theta \in \Theta}
\det C_0(\theta) >0$ and {\bf (Id)} holds, then the strong consistency holds.\\
Moreover if $\theta_0\in  \stackrel{\circ}{ \Theta}$ with
$\stackrel{\circ}{ \Theta}\subset \widetilde\Theta(4)$, then
$\widehat\theta_n$ is also asymptotically normal.
\end{Proposition}
The asymptotic normality was settled before by Comte and Lieberman
\cite{Comte2003} for $r=8$. Our result needs just moment of order
$r=4$.
\subsection{Multivariate {\em NLARCH($\infty$)} models}
Let $\big ( B_j(\theta) \big )_{j \ge 1}$
be a sequence of $m \times d$-matrices  and $B_0(\theta)$ be a vector of $\R^m$. The multivariate LARCH($\infty$) models
introduced by
Doukhan {\it et al.} \cite{Doukhan2006} are extensions of the univariate LARCH($\infty$) models of Giraits {\it et al.} \cite{Giraitis2004}. They are defined as the stationary solution
of the equation:
\begin{equation*}
X_t = \zeta_t \Big( B_0(\theta_0) + \underset{j =
1}{\overset{\infty}{\sum}}B_j(\theta_0) X_{t - j} \Big).
\end{equation*}
Notice that the innovations $(\zeta_t)_{t\in\Z}$ are here random matrices. In this context, the QMLE is not a
suitable estimator since $\inf_{\theta\in\Theta}\det
\big(H_\theta(x) \big ) =0$ except in very
specific cases. However, Doukhan and
Wintenberger \cite{Doukhand} proposed a generalization of
LARCH$(\infty)$ models, so-called NLARCH$(\infty)$ models defined by the equation
\begin{eqnarray}
\label{e:arch_infty}X_t&=&\zeta_t\Big(
B_0(\theta)+\sum_{j=1}^\infty B_j(\theta,X_{t-j})\Big),
\end{eqnarray}
where now $B_j(\theta,.):\R^m\to\R^p$ are $b_j(\theta)$-Lipschitz
functions. If the matrices of the innovations are
concentrated on the diagonal, we rewrite \eqref{e:arch_infty} as
$$
X_t = M_{\theta_0}(X_{t-1},X_{t-2},\ldots)\xi_t,
$$
where $(\xi_t)_i=(\zeta_t)_{i,i}$ and $\big (
M_{\theta}(X_{t-1},X_{t-2},\ldots)\big )_{ij}=\delta_{ij} \cdot \big (
B_0(\theta) + \underset{k = 1}{\overset{\infty}{\sum}} B_k(\theta,
X_{t - k}) \big)_i$. For
instance, consider the multidimensional extension of the TARCH models as
$$(B_j(\theta,x))_k=\sum_{i=1}^mB^+_{j,k,i}(\theta)
\max(x_{j,i},0)+B^-_{j,k,i}(\theta)\min(x_{j,i},0),
$$
where $B^+_{j,k,i}$ and $B^-_{j,k,i}$ are non-negative real
numbers. For NLARCH($\infty$) models we have
$$
\Theta(r)=\left\{\theta\in\R^d~\Big/~\sum_{j=1}^{\infty}\Lip_x
(B_j(\theta,x))<\Big (
\E \big [ \| \xi_0\|^r \big ] \Big )^{-1/r}\right\}.
$$
\begin{Proposition}
Let $\Theta$ be a compact subset of $\Theta(2)$, $X$ be the stationary solution to
\eqref{e:arch_infty} and assume that {\bf (Id)} holds.
Assume that for all $j \in \N$, the vectors
$B_j(\theta,.)\in [0,\infty[^p$, $ \inf_{\theta \in \Theta}
\|B_0(\theta) \| >0$ and for all $j \in \N^*$,
$$
\|\Lip_x (B_j(\theta,x))\|_\Theta =O\big
(j^{-\ell}\big )~~\mbox{for some}~~\ell >3/2,
$$
then the strong consistency holds.\\
Moreover if $\theta_0\in  \stackrel{\circ}{ \Theta}$ with
$\stackrel{\circ}{ \Theta}\subset \widetilde\Theta(4)$, assume that
the functions $\theta \mapsto B_j(\theta,.)$ are $2$ times
continuously differentiable on $\Theta$ and satisfy
\begin{multline*}
\Big \|\lip_x \frac {\partial
B_j(\theta,x)} {\partial \theta_k } \Big \|_\theta =O\big (j^{-\ell'}\big
)~~\mbox{with}~~\ell'>3/2\\
\mbox{ and }\sum_{j\ge1}\Big \|\lip_x \frac {\partial^2
B_j(\theta,x)} {\partial \theta_k
\partial \theta_{k'}} \Big \|_\Theta<\infty\mbox{  for all }(k,k')\in
\{1,\ldots,d\}^2.
\end{multline*}
If {\bf (Var)} holds, the QMLE $\widehat\theta_n$ is
also asymptotically normal.
\end{Proposition}

\subsection{Multivariate non-linear {\em AR($\infty$)} models}
Let us focus on one example where $M=I_m$ and $f\neq 0$. In this
context, {\bf (A1(M))} is always satisfied and the QMLE coincides with the
least squares error estimator. Here, we restrict ourselves to the cases where
\begin{equation}\label{eq::ar}
f_\theta(X_{t-1},X_{t-2},\ldots)=A_0(\theta_0)+\sum_{i=1}^\infty
A_{i}(\theta_0,X_{t-i}),
\end{equation}
where $A_{i}(\theta,.)$ are Lipschitz functions with values in
positive definite $d \times d$ matrices. Here $\Theta(r)$ neither
depends on $r$ nor on the distribution of $\xi_0$:
$$
\Theta(r)=\left\{\theta\in\R^d~\Big/~\sum_{j=1}^\infty  \lip_x
A_j(\theta,x)<1\right\}. $$
\begin{Proposition}
Let $\Theta$ be a compact subset of $\Theta(2)$, $X$ be the stationary solution
to \eqref{eq::sys} when $M=I_m$ and \eqref{eq::ar} holds. Under {\bf (Id)} and if for $j \in
\N^*$,
$$
\left\|\lip_x A_j\right\|_\Theta=O\big (j^{-\ell}\big
)~~\mbox{for some}~~\ell>3/2,
$$
then the strong consistency
holds.\\
Moreover if $\theta_0\in  \stackrel{\circ}{ \Theta}$ with
$\stackrel{\circ}{ \Theta}\subset \widetilde\Theta(4)$,
$\widehat\theta_n$ is also asymptotically normal as soon as functions $\theta \mapsto A_j(\theta,x)$ are
$2$ times continuously differentiable on $\Theta$ for all
$x\in\R^m$, {\bf (Var)} holds and
\begin{multline*}
\Big\|\lip_x\frac{\partial A_j}{\partial \theta_k}\Big
\|_\Theta=O\big (j^{-\ell'}\big
)~~\mbox{for some}~~\ell'>3/2\mbox{ and}\\
\sum_{j=1}^\infty \Big \|\lip_x\frac{\partial^2 A_j}{\partial
\theta_k\theta_{k'}}\Big\|_\Theta<\infty\mbox{ for all
}k,k'\in\{1,\ldots,d\}.
\end{multline*}
\end{Proposition}

\subsection{Multivariate {\em ARMA-GARCH} models}\label{armagarch}
Here $M_\theta$ is concentrated on its diagonal and $f$ is not necessarily
identically zero. If $f\equiv0$, the model coincides with the VEC-GARCH
model, see Jeantheau
\cite{Jeantheau1998}. Multidimensional ARMA-GARCH processes were
introduced by Ling and McAleer \cite{Ling2003} as the solution of the system of equations
\begin{equation}\label{arma}
\begin{cases}
&\Phi_\theta(L) \cdot X_t=\Psi_\theta(L) \cdot \varepsilon_t,\\
&\varepsilon_t=M_\theta(X_{t-1},X_{t-2},\ldots)\xi_t,
\end{cases}
\end{equation}
with diag$(H_\theta^t)=C_0(\theta)+\sum_{i=1}^qC_i(\theta)$diag$(\varepsilon_{t-i}\varepsilon_{t-i}')+
\sum_{i=1}^{q'}D_i(\theta)$diag$(H_\theta^{t-i})$.
Here $C_0(\theta)$, $C_i(\theta)$ and $D_j(\theta)$ are positive definite matrices,
$\mbox{diag}A$ is the diagonal of the matrix $A$,
$\Phi_\theta(L)=I_m-\Phi_1L-\cdots-\Phi_sL^s$ and
$\Psi_\theta(L)=I_m-\Psi_1L-\cdots-\Psi_{s'}L^{s'}$ are polynomials
in the lag operator $L$ and $\Phi_i$ and $\Psi_j$ are squared
matrix.
We define for all $\theta\in\R^d$
\begin{eqnarray*}
\Gamma_\theta(L)&:=&I_m+\sum_{i=1}^\infty\Gamma_i(\theta)L^i=\Psi_\theta^{-1}(L)\Phi_\theta(L) \\
\mbox{and}~~\sum_{i=1}^\infty B_i(\theta)  Z^i &:=& \Big(1 -
\sum_{i=1}^{q'} D_i(\theta) Z^i\Big)^{-1}\times\sum_{i=1}^q C_i(\theta)  Z^i~~\mbox{for
all}~~Z\in \C^m,
\end{eqnarray*}
where the polynomials of the right hand side are assumed to be
coprime. The equation (\ref{arma}) has the representation \eqref{eq::sys}
with $f_\theta(X_{t-1},X_{t-2},\ldots):=\sum_{i=1}^\infty\Gamma_i(\theta)
X_{t-i}$. We can define,
$$
\Theta(r)=\left\{\theta\in\R^d~\Big/~
\sum_{i=1}^\infty\|\Gamma_i(\theta)\|+\Big ( \E \big [ \| \xi_0\|^r \big
] \Big )^{1/r}\sum_{j=1}^\infty \|B_j(\theta)\|< 1\right\}.
$$
If $\theta_0\in\Theta(r)$ then the existence of a solution is
ensured. This existence condition is more explicit than the one of
Theorem 2.1. of Ling and McAleer \cite{Ling2003}. Now we give a version of Theorems \ref{th::consist} and
\ref{th::asnorm} when
$$
\theta=(\Phi_1,\ldots,\Phi_s,\Psi_1,\ldots,\Psi_{s'},
C_0,C_1,\ldots,C_q,D_1,\ldots,D_{q'}).
$$
\begin{Proposition}\label{as::armagarch}
Let $\Theta$ be a compact subset of $\Theta(2)$, $X$ be the stationary solution
to the system \eqref{arma}. If $\inf_{\theta \in \Theta} \det
C_0(\theta) >0$ and {\bf (Id)} holds then $\widehat\theta_n$ is
consistent.\\
Moreover if $\theta_0\in  \stackrel{\circ}{ \Theta}$ with
$\stackrel{\circ}{ \Theta}\subset \widetilde\Theta(4)$,
$\widehat\theta_n$ is also asymptotically normal as soon as {\bf (Var)} holds.
\end{Proposition}
Ling and McAleer \cite{Ling2003} also provided consistency and
asymptotic normality of the QMLE.
Proposition \ref{as::armagarch} improves their results. Notice that for
VEC-GARCH models, Jeantheau \cite{Jeantheau1998} provided the consistency under a weaker
condition.

\section{Proofs} \label{PROOFS}
In this section the proofs of the main results are collected in
the order of appearance in the paper. First we prove Proposition
\ref{stationarity} and Corollary \ref{scal}, then Lemma
\ref{lem::estvar} that settles the invertibility of the QMLE. With the help of this property we
prove the main theorems that state consistency and asymptotic
normality of the QMLE.
\subsection{Proofs of Proposition \ref{stationarity} and Corollary \ref{scal}}\label{proof1}
We apply a result of Douk\-han and Wintenberger \cite{Doukhand} that gives conditions
for the existence of a stationary solution of an equation of type
\begin{eqnarray} \label{pdow}
X_t=F\big (X_{t-1},X_{t-2},\ldots;\xi_t \big )~~\mbox{a.e. for all}~~t\in
\Z.\end{eqnarray}If $\E \|\xi_0\|^r< \infty$ and
$F$ satisfies for $x=(x_i)_{i \ge1}$, $y=(y_i)_{i \ge1}\in
({\R^m}) ^{\infty}$,
\begin{itemize}
\item $\E\|F(0;\xi_0)\|^r <
\infty$;
\item $\displaystyle \Big ( \E\|F\big (x\ ; \xi_0 \big )-F\big (y\ ;  \xi_0 \big)\|^r
\Big ) ^{1/r}  \hspace{-2mm} \leq \sum_{j \ge1} a_j
\|x_j-y_j\|$, with $\displaystyle \sum_{j\ge1}a_j <1$;
\end{itemize}
the existence of a unique causal stationary solution $X$ of (\ref{pdow}), such that
$\E\big [\|X_0\|^r\big ]<\infty$ is proved in \cite{Doukhand}. We identify $F$ from (\ref{eq::sys}):
$$
F\big (X_{t-1},X_{t-2},\ldots;\xi_t \big )=
M_{\theta_0}(X_{t-1},X_{t-2},\ldots) \cdot \xi_t+
f_{\theta_0}(X_{t-1},X_{t-2},\ldots).
$$
Obviously, $\E\big [
\|F(0;\xi_0)\|^r \big ] < \infty$ if $\E\|\xi_0\|^r<\infty$ and we have
\begin{eqnarray*}
\Big ( \E\big \|F\big (x\ ; \xi_0 \big )-F\big
(y\ ;  \xi_0 \big )\big\|^r \Big )
^{1/r}&&\\
&&\hspace{-2.5cm} \leq\Big ( \E\big
\|(M_{\theta_0}(x)- M_{\theta_0}(y))\cdot \xi_0 \big
\|^r \Big )
^{1/r}+  \big \| f_{\theta_0}(x)-f_{\theta_0}(y)  \big \| \\
&& \hspace{-2.5cm} \leq \big(\E \|\xi_0\|^r\big
)^{1/r}\big \| M_{\theta_0}(x)-M_{\theta_0}(y)
\big \|+\big \| f_{\theta_0}(x)-f_{\theta_0}(y)
\big \|.
\end{eqnarray*}
The condition of Proposition \ref{stationarity} then implies those
of \cite{Doukhand} on $F$. In the context of Corollary \ref{scal},
from $H_{\theta}(x)=\widetilde H_{\theta}(x^2)$ for all
$x=(x_j)_{j\ge1}\in \R^\infty$, we have
\begin{eqnarray*}\|\widetilde H_{\theta}(x)-\widetilde
H_{\theta}(y)\|
&\le& \sum_{j=1}^\infty \alpha_j(H,\theta)\|x_j-y_j\|\qquad \mbox{ and,}\\
\big(\E\big [
|M_{\theta_0}^2(x)\xi_0^2-M_{\theta_0}^2(y)\xi_0^2|^{r/2}\big
]\big)^{2/r}&= &\big(\E\big [|\xi_0|^r\big ]\big)^{2/r}|\widetilde
H_{\theta_0}(x^2)
-\widetilde H_{\theta_0}(y^2)|\\
&\le&\big(\E\big [|\xi_0|^r\big ]\big)^{2/r}\sum_{j=1}^\infty
\alpha_j^0(H,\theta_0)|x_j^2-y_j^2|.
\end{eqnarray*}
The results of \cite{Doukhand} yield the existence in $\L^{r/2}$ of
the solution $(X_t^2)_{t\in\Z}$ of the equation
$$
X_t^2=M_{\theta_0}^2(X_{t-1},X_{t-2},\ldots)\xi_t^2= \widetilde
H_{\theta_0}(X^2_{t-1},X^2_{t-2},\ldots)\xi_t^2\qquad
\mbox{a.e.}
$$
Moreover, by \cite{Doukhand} there exists a measurable function
$\varphi$ such that $X_t=\varphi(\xi_t,\xi_{t-1},\ldots)$ for all
$t\in\Z$. The ergodicity of $X$ follows from the Proposition 4.3 in
Krengel \cite{Krengel1985}; it states that if $(E,{\cal E})$ and
$(\widetilde E,\widetilde {\cal E})$ are measurable spaces,
$(v_t)_{t\in\Z}$ is a stationary ergodic sequence of $E$-valued
random elements and $\varphi:(E^{\N}, {\cal E}^{\N})\mapsto
(\widetilde E,\widetilde {\cal E})$ is a measurable function then
the sequence $(\widetilde v_t)_{t\in\Z}$ defined by $\widetilde
v_t=\varphi(v_t,v_{t-1},\ldots)$ is a stationary ergodic process.

\subsection{Proof of Lemma \ref{lem::estvar}}
We treat the three assertions of the lemma one after the other.\\
1. Define $f_\theta^{t,p}=f_\theta(X_{t-1},\ldots,X_{t-p},0,0,\ldots)$
for all $t\in\Z$ and $p\in \N$. We have
$f_\theta^{t,p}\in \L^r({\cal C}(\Theta,\R^m))$ because $\theta_0\in\Theta(r)$ and, using Corollary \ref{stationarity}, all the following quantities are finite:
\begin{eqnarray*}
\Big(\E \big [\|f_\theta^{t,p}\|_\Theta ^r \big ] \Big)^{1/r} &\leq & \Big(\E \big
[\|f_\theta^{t,0}-f_\theta^{t,p}\|_\Theta^r \big ]\Big)^{1/r} + \Big(\E \big
[\|f_\theta^{t,0}\|_\Theta^r \big ]\Big)^{1/r} \\
& \leq & \Big (\sum_{j\ge
1}\alpha_j(f)\Big )\Big(\E\big [\|X_0\|^r\big ]\Big)^{1/r}  + \|f_\theta(0)\|_\Theta.
\end{eqnarray*}
For $p<q$
$$
\E \big [\|f_\theta^{t,p}-f_\theta^{t,q}\|_\Theta^r \big ] \leq
\E \Big [ \Big \|\sum_{p< j\le q} \alpha_j(f) X_{t-j} \Big \| ^r
\Big ]\leq \E\big [ \|X_0\|^r\big ] \Big (\sum_{p< j\le q}\alpha_j(f)\Big )^r.
$$
Since $\sum_{j\ge 1}\alpha_j(f)<\infty$, $(f_\theta^{t,p})_{p\ge
0}$ satisfies the Cauchy criteria in $\L^r({\cal C}(\Theta,\R^m))$
and it converges to $f_\theta^{t,\infty}$, that is $f_\theta^t$ on
$\sigma(X_{t_1},\ldots,X_{t_n})$ for all $n\in\N^\ast$ and
$t>t_1>\cdots>t_n$ (those $\sigma$-algebras generate
$\sigma(X_{t-1},X_{t-2},\ldots)$
and therefore $f_\theta^{t,\infty}=_{a.s}f_\theta^t$). \\

2. Define $H^{t,p}_\theta=
H_\theta(X_{t-1},\ldots,X_{t-p},0,\ldots) $ for all $p\in \N$ and
$t\in \N$. From Corollary \ref{stationarity}, $\theta_0\in\Theta(r)$ and common inequalities
satisfied by matrix norms,  $H^{t,p}_\theta\in\L^{r/2}({\cal
C}(\Theta,{\cal M}_m))$ since, denoting
$M_\theta^{t,p}=M_\theta(X_{t-1},\ldots,X_{t-p},0,\ldots)$,
$$
\|H^{t,p}_\theta\|^{r/2}_\Theta \le
\|M_\theta^{t,p}\|_\Theta^r \le\Big (\|M_\theta(0)\|_\Theta+\sum_{j=1}^{\infty}
\|X_{t-j}\|\alpha_j(M)\Big)^r.
$$
We conclude as above that $H_\theta^t\in \L^{r/2}({\cal
C}(\Theta,{\cal M}_m))$ by bounding, for $p<q$,
\begin{eqnarray*}
\|H^{t,p}_\theta-H^{t,q}_\theta\|^{r/2}_\Theta &\hspace{-3mm} \le&\hspace{-3mm}
\big \|M_\theta^{t,p}
-M_\theta^{t,q}\big \|_\Theta^{r/2}\Big ( \|
M_\theta^{t,p}\|^{r/2}_{\Theta
}\hspace{-2mm}+ \hspace{-1mm} \|
M_\theta^{t,q}\|^{r/2}_{\Theta } \Big
).
\end{eqnarray*}
The Cauchy-Schwarz inequality implies that
\begin{eqnarray*}
\E \big [\|H_\theta^{t,p}-H_\theta^{t,q}\|^{r/2}_\Theta \big ]&&\hspace{-6mm}\le\Big ( \E \big [
\|(M_\theta^{t,p}-
M_\theta^{t,q}\big \|_\Theta^ {r} \big
] \Big )^{1/2}\\
&&\hspace{-6mm}\times\Big[\Big (\E \big [
\|M_\theta^{t,p}\|_\Theta^{r}\big ]\Big )^{1/2} +
\Big(\E \big [
\|M_\theta^{t,q}\|_\Theta^{r}\big ]
\Big )^{1/2}\Big] \\
&&\hspace{-6mm} \le B\Big( \E \Big [ \Big
(\sum_{p< j\le q}\alpha_j(M)
\|X_{t-j}\|\Big)^r \Big ] \Big)^{1/2}\\
&&\hspace{-6mm}\le B \left(\E \big [ \|X_{0}\|^r
\big ] \right)^{1/2}\Big ( \sum_{p< j\le q}\alpha_j(M)\Big
)^{r/2}
\end{eqnarray*}
for some constant $B>0$. \\

3. First notice that $\|X_0X'_0\|\le\|X_0\|^2$. Next, as in the
previous proofs, $(H^{t,p}_\theta)_{p\in\N^\ast}$ converges to
$H_\theta^t$ in $\L^{r/2}({\cal C}(\Theta,{\cal M}_m))$. Thus there exists a subsequence $(p_k)_{k\in\N}$
such that $\| H_\theta^{t,p_k}- H^t_\theta \|_\Theta  \limitepsk
0$. Thanks to the continuity of the determinant, $(\det
H_\theta^{t,p_k})_{k\in\N}$ also converges a.s. to $\det
H_\theta^{t}$. Then $\det H_\theta^t\ge \underline H$,
$H_\theta^t$ is an invertible matrix and in view of elementary relations
between matrix norm and determinant $\Big\|\big (\widehat
H_\theta^t  \big) ^{-1}\Big\|_\Theta\le \underline H^{-1/m}$.

\subsection{Proof of Theorem \ref{th::consist}}\label{proof2}
The proof of the theorem is divided into two parts. In (i) a uniform
(in $\theta$) law of large numbers on $(\widehat q_t)_{t \in \N^*}$
(defined in (\ref{eq::QMLE})) is established. In (ii), it is
proved  that $L(\theta):=-\E(q_t(\theta))/2$ has a unique maximum
in $\theta_0$. Those two conditions lead to
the consistency of $\widehat \theta_n$.\\

(i) Using Proposition \ref{stationarity}, with
$q_t=G(X_{t},X_{t-1},\cdots)$, one deduces that $(q_t)_{t\in\Z}$
(defined in (\ref{eq::MLE})) is a stationary ergodic sequence.
>From Straumann and Mikosch \cite{Straumann2006}, we know that if
$(v_t)_{t\in\Z}$ is a stationary ergodic sequence of random
elements with values in $\C(\Theta,\R^m)$, then the uniform (in
$\theta \in \Theta$) law of large numbers is implied by
$\E\|v_0\|_\Theta<\infty$. As a consequence, $(q_t)_{t\in\Z}$
satisfies a uniform (in $\theta \in \Theta$) strong law of large
numbers as soon as $\E \big [ \sup_\theta|q_t(\theta)|\big ]
<\infty$. But, from the inequality $\log(x) \leq x-1$ for all $x
\in ]0,\infty[$ and Lemma \ref{lem::estvar}, for all $t\in \Z$,
\begin{eqnarray}
\nonumber &&|q_t(\theta)|\le \frac
{\|X_t-f_t(\theta)\|^2}{(\underline H)^{1/m}}+m \Big| \frac 1 m
\log \underline H+ \frac {\|H_\theta^t \|}{{\underline
M}^{1/m}}-1\Big| ~~\mbox{for all $\theta \in \Theta$} \\
&&\hspace{1cm}\label{qborne} ~~\Longrightarrow~~\sup_{\theta \in
\Theta} |q_t(\theta)| \leq  \frac
{\|X_t-f_t(\theta)\|_\Theta^2}{(\underline H)^{1/m}}+\Big | \log
\underline H \Big | + m   \times \frac {\|H_\theta^t
\|_\Theta}{{\underline H}^{1/m}}.
\end{eqnarray}
But for all $t\in \Z$, $\E \|X_t\|^r<\infty$, see Corollary \ref{stationarity}, and $\E\big [
\|f^t_\theta\|^r_\Theta\big ]+\E \big [\|H_\theta^t \|_\Theta^{r/2} \big] <\infty$, see
Lemma \ref{lem::estvar}. As a consequence, the right hand side of
(\ref{qborne}) has a finite first moment
and therefore
$$
\E \big [ \sup_{\theta\in\Theta}|q_t(\theta)| \big ]<\infty.
$$
The uniform strong law of large numbers for $(q_t(\theta))$ directly
follows and hence
\begin{eqnarray}\label{convLn}
\Big \|\frac { L_n(\theta)} n-L(\theta) \Big \|_\Theta \limiteasn
0~~~\mbox{with}~~L(\theta):=-\frac 1 2 \E \big [q_0(\theta)\big ].
\end{eqnarray}
Now, one shows that $\displaystyle \frac 1 n \|\widehat
L_n-L_n\|_\Theta \limiteasn 0$. Indeed, for all $\theta \in
\Theta$ and $t\in \N^*$,
\begin{eqnarray}
\nonumber  \lefteqn{\big |\widehat q_t(\theta)-q_t(\theta)\big|}\\
\nonumber &= & \Big ( \log \det \widehat H_\theta^t - \log
\det H_\theta^t \Big )+(X_t-\widehat f^t_\theta\big)'
\big ( \widehat H_\theta^t  \big) ^{-1} \big (X_t-\widehat
f^t_\theta\big)\\
\nonumber &&- (X_t-f^t_\theta\big)'\big ( H_\theta^t \big)
^{-1} \big (X_t-f^t_\theta\big)\\
\nonumber &&  \leq|C|^{-1} \Big | \det \big(\widehat
H_\theta^t \big )- \det
\big(H_\theta^t \big ) \Big |+(X_t-\widehat
f^t_\theta\big)'\Big [ \big ( \widehat H_\theta^t  \big) ^{-1}
\hspace{-3mm} -\hspace{-1mm} \big ( H_\theta^t  \big) ^{-1}  \Big
]\big (X_t-\widehat f^t_\theta\big)\\
\nonumber && +(2X_t-\widehat f^t_\theta-
f^t_\theta \big)'\big ( H_\theta^t  \big) ^{-1}  \big
(f^t_\theta-\widehat
f^t_\theta\big)\\
\nonumber &\leq&
{\underline H}^{-1}\Big \| \det \big(\widehat H_\theta^t \big
)- \det \big(H_\theta^t \big ) \Big \|_\Theta +2 \big (\|X_t \|+\|\widehat
f^t_\theta \|_\Theta \big )\big \|\big ( \widehat H_\theta^t \big)
^{-1}\hspace{-3mm} -\big ( H_\theta^t \big) ^{-1} \big \|_\Theta \\
\label{inegq} &&  + \Big (2 \|X_t\|+\|\widehat
f^t_\theta\|_\Theta +\|f^t_\theta\|_\Theta \Big )\big \| \big (
H_\theta^t  \big) ^{-1} \big \|_\Theta \big \|f^t_\theta-\widehat
f^t_\theta\big\|_\Theta
\end{eqnarray}
% because when $A$ is definite positive matrix $\|x \|_A^2:=x'  A x$ is a norm such that $|\|x_1\|_A^2-\|x_2\|_A^2 | \leq \|x_1-x_2 \|^2_A$.
by the mean value theorem, with $C \in [ \det \big(H_\theta^t
\big ),\det \big(\widehat H_\theta^t \big )]$ and therefore
$|C|>\underline H$. On the one hand,
\begin{eqnarray*}
\|\big ( \widehat H_\theta^t  \big) ^{-1}-\big ( H_\theta^t  \big)
^{-1} \big \|_\Theta &\leq  & \big \|\big ( \widehat H_\theta^t
\big) ^{-1} \big \|_\Theta  \big \| \widehat H_\theta^t
-H_\theta^t  \big \|_\Theta \cdot  \big \| \big ( H_\theta^t  \big)
^{-1} \big \|_\Theta.
\end{eqnarray*}
On the other hand, for  an invertible
matrix $A \in {\cal M}_m(\R)$, and $H  \in {\cal M}_m(\R)$,
$$
\det (A+H) = \det (A)+\det(A) \cdot \mbox{Tr}\big ((A^{-1})' H\big
)+o(\|H\|),
$$
where $ \big | \mbox{Tr}\big ((A^{-1})' H\big )\big | \leq \big
\|A^{-1}\big \| \cdot \big \| H \big \|$. Using the relation $\| \big ( H_\theta^t  \big) ^{-1} \big
\|_\Theta \geq \underline H^{-m}$ for all $t\in\Z$, there exists $C>0$ not
depending on $t$ such that inequality (\ref{inegq}) becomes:
$$
\sup_{\theta \in \Theta} \big |\widehat
q_t(\theta)-q_t(\theta)\big|\leq C
\big (\|X_t \|+\|\widehat f^t_\theta \|_\Theta +\|
f^t_\theta\|_\Theta \big ) \times \Big ( \big \| \widehat H_\theta^t
-H_\theta^t  \big \|_\Theta+ \big \|f^t_\theta-\widehat
f^t_\theta\big\|_\Theta \Big )
$$
>From the H\"older and Minkowski inequalities and by virtue of $3/2=
1+1/2$,
\begin{eqnarray}
\nonumber \E \big [\sup_{\theta \in \Theta} \big |\widehat
q_t(\theta)-q_t(\theta)\big |^{2/3} \big ] &&\hspace{-6mm} \leq
 C\Big(\E  \big [\|X_t \|+\|\widehat
f^t_\theta \|_\Theta +\| f^t_\theta\|_\Theta \big]^{2}\Big)^{1/3}\\
\nonumber&&\hspace{-6mm}\times\Big( \E\big [\| \widehat H_\theta^t -H_\theta^t
\|_\Theta\big
]+ \E\big [\|f^t_\theta-\widehat
f^t_\theta\|_\Theta\big
]\Big)^{2/3}\\
&& \hspace{-3mm}\label{borne2/3}\hspace{-3mm}  \leq  C'\Big ( \sum_{j\ge t}
\big [ \alpha_j(f)+\alpha_j(M) \big ]\Big )^{2/3},
\end{eqnarray}
with $C'>0$ not
depending on $\theta$ and $t$. Now,
consider for $n \in \N^*$,
$$
S_n:=\sum_{t=1}^n\frac{1}{t}\sup_{\theta \in \Theta} \big
| \,\widehat q_t(\theta)-q_t(\theta)\big |.
$$
Applying the Kronecker lemma (see Feller \cite{Feller2}, p. 238),
if $ \lim_{n \to \infty} S_n<\infty$ a.s.  then $\displaystyle
\frac 1 {n} \cdot \|\widehat L_n-L_n\|_\Theta \limiteasn 0$. Following
Feller's arguments, it remains to show that for all
$\varepsilon>0$,
\begin{eqnarray*}
\label{probaA} \P(\forall n\in \N,~~\exists m>n\mbox{ such that
}|S_m-S_n|>\varepsilon):= \P(A)=0.
\end{eqnarray*}
Let $\varepsilon>0$ and denote
$$
A_{m,n}:=\{|S_m-S_n|>\varepsilon\}
$$
for $m>n$. Notice that $A= \bigcap_{n\in\N}\bigcup_{m>n}A_{m,n}$.
For $n \in \N^*$, the sequence of sets $(A_{m,n})_{m>n}$ is
obviously increasing, and if $A_n:=\bigcup_{m>n}A_{m,n}$, then
$\lim_{m\to\infty}\P(A_{m,n})=\P(A_n)$. Observe that $(A_n)_{n\in\N}$
is a decreasing sequence of sets and thus,
$$
\lim_{n\to\infty}
\lim_{m\to\infty}\P(A_{m,n})=\lim_{n\to\infty}\P(A_n) =\P(A).
$$
It remains to bound $\P(A_{m,n})$. From the Bienaym\'e-Chebyshev
inequality,
\begin{eqnarray*}
\P(A_{m,n})&&\hspace{-6mm}=\P\Big(\sum_{t=n+1}^m \frac{1}{t} \sup_{\theta \in
\Theta} \big |\widehat q_t(\theta)-q_t(\theta)\big
|>\varepsilon\Big)\\&&\hspace{-6mm}\le \frac{1}{\varepsilon^{2/3}}\E\Big[\Big
(\sum_{t=n+1}^m \frac{1}{t^{}}\sup_{\theta \in \Theta} \big
|\widehat
q_t(\theta)-q_t(\theta)\big | \Big )^{2/3}\Big]\\
&&\hspace{-6mm}\le \frac{1}{\varepsilon^{2/3}}\sum_{t=n+1}^m\frac{1}{t^{2/3}}\E
\Big [ \sup_{\theta \in \Theta} \big |\widehat
q_t(\theta)-q_t(\theta)\big |^{2/3}\Big ].
\end{eqnarray*}
Using (\ref{borne2/3}) and condition (\ref{condP}), since
$\ell>3/2$, there exists $C>0$ such that $$\displaystyle \Big (\sum
_{j=t}^\infty \alpha_j(f)+\alpha_j(M)+\alpha_j(H) \Big ) ^{2/3} \hspace{-3mm}
\leq \hspace{-1mm} \frac C {t^{2(\ell-1)/3}}.$$
Thus,
$\displaystyle t^{-2/3}\E \Big [\sup_{\theta \in \Theta}
\big |\widehat
q_t(\theta)\hspace{-1mm}-\hspace{-1mm}q_t(\theta)\big |^{2/3}\Big
]\hspace{-1.5mm}\le C\big(t^{-2\ell/3}\big)$ for some $C>0$ and
$$
\sum_{t=1}^\infty\frac{1}{t^{2/3}}\E \Big [
\sup_{\theta \in \Theta} \big |\widehat
q_t(\theta)-q_t(\theta)\big |^{2/3} \Big ] <\infty\mbox{ as } \ell>3/2.
$$
Thus $\lim_{n\to\infty}
\lim_{m\to\infty}\P(A_{m,n})\limiten 0$
and $\displaystyle \frac 1 {n} \cdot \|\widehat
L_n-L_n\|_\Theta \limiteasn 0$.\\

(ii) See Proposition 2.1. of Jantheau \cite{Jeantheau1998}.

\subsection{Proof of Theorem \ref{th::asnorm}}\label{proof3}
Let $V$ be a Banach space (thereafter $V=\R^m$ or $V={\cal M}_m$)
and ${\cal D}^{(2)}{\cal C}(\Theta,V)$ denote the Banch space of
$V$-valued $2$ times continuously differentiable functions on
$\Theta$ equipped with the uniform norm
$$
\displaystyle \|g\|_{(2),\Theta}=\|g\|_\Theta+ \Big
\|\frac {\partial g}{\partial \theta} \Big
\|_\Theta+\Big \|\frac {\partial^2
g}{\partial \theta \partial \theta'} \Big \|_\Theta.
$$
We start by proving the following preliminary lemma:
\begin{Lemma}\label{lem::prel}
Let $\theta_0$ belong to $\Theta(r)$ ($r\ge2$) and assume that {\bf (A3(f))} and {\bf (A3(M))} or {\bf (A3(H))} hold. Then $$f_\theta^t\in\L^{r}\big ({\cal D}^{(2)}{\cal
C}(\Theta,\R^m)\big )\mbox{ and }H_\theta^t\in\L^{r/2}\big ({\cal D}^{(2)}{\cal
C}(\Theta,{\cal M}_m)\big ).$$
\end{Lemma}
In view of the results of Lemmas \ref{lem::estvar} and
\ref{lem::prel}, the functions ${\partial
L_n(\theta)}/{\partial\theta}$ and ${\partial^2
L_n(\theta)}/{\partial\theta^2}$ are measurable and a.s. finite
for all $\theta\in {\Theta}$. Their asymptotic properties are described in the next two lemmas
\begin{Lemma}\label{lem::mart}
Let $\theta_0$ belong to $\Theta(r)$ ($r\ge4$) and assume that {\bf (A3(f))} and {\bf (A3(M))} or {\bf (A3(H))} hold, then
\begin{equation}\label{clt::stat}
n^{-1/2} \frac {\partial L_n(\theta_0)}{\partial \theta}
\limiteloin  {\cal N}_d(0,G(\theta_0)),
\end{equation}
where $G(\theta_0)=(G(\theta_0))_{1\leq i,j \leq d}$ is finite and its expression is
given in \eqref{G0}.
\end{Lemma}
\begin{Lemma}\label{lem::der2}
Let $\theta_0$ belong to $\Theta(r)$ ($r\ge4$) and assume that {\bf (A3(f))} and {\bf (A3(M))} or {\bf (A3(H))} hold, then
\begin{equation}\label{convD2Ln}
\Big \|\frac 1n\frac{\partial^2 L_n(\theta)}{\partial \theta\partial \theta'}
-\frac{\partial^2 L(\theta)}{\partial \theta\partial \theta'} \Big \|_\Theta \limiteasn
0~~~\mbox{with}~~\frac{\partial^2 L(\theta)}{\partial \theta\partial \theta'}:=
-\frac 1 2 \E \Big [\frac{\partial^2 q_0}{\partial \theta\partial \theta'}(\theta)\Big ].
\end{equation}
\end{Lemma}
We postponed the proofs of Lemmas \ref{lem::estvar}-\ref{lem::der2} to the end of the Section and continue with the proof of Theorem \ref{th::asnorm}. From  Theorem
\ref{th::consist}, we have
\begin{eqnarray}\label{conv:thetan}
\widehat \theta_n\limiteasn \theta_0.
\end{eqnarray}
Since $\theta_0\in \stackrel{\circ}{\Theta}$, a Taylor expansion
of ${\partial L_n(\theta_0)}/{\partial \theta_i}\in\R$ implies
\begin{equation}\label{eq::tayl}
\frac {\partial L_n(\widehat \theta_n)}{\partial \theta_i} = \frac
{\partial L_n( \theta_0)}{\partial \theta_i}+\frac {\partial^2
L_n(\overline \theta_{n,i})}{\partial \theta\partial \theta_i}(\widehat
\theta_n-\theta_0),
\end{equation}
for $n$ sufficiently large such that the $\overline\theta_{n,i}\in
\Theta$, which are between $\widehat \theta_n$ and $\theta_0$ for
all $1\le i\le d$. Using equations (\ref{convD2Ln}) and
(\ref{conv:thetan}), we conclude with the uniform convergence
theorem that
$$
F_n:=-2\Big(\frac 1 n\frac {\partial^2
L_n(\overline \theta_{n,i})}{\partial \theta\partial \theta_i}\Big)_{1\le i\le d}\limiteasn
F(\theta_0).
$$
One obtains $\big (F(\theta_0)\big )_{ij} =\E\Big [{\partial^2
q_0(\theta_0)}/{\partial \theta_i
\partial \theta_j }\Big]$ for $1\le i,j\le d$. With similar arguments
as for (\ref{calcder}), since $X_t-f_{\theta_0}^t
=M_{\theta_0} \xi_t$, with $\xi_t$ independent of
$(X_{t-1},X_{t-2},\ldots)$,
\begin{multline*}
\E\Big[\big (
X_t-f_{\theta_0}^t \big)' \frac {\partial^2
\big(H_\theta^t\big)^{-1}}{\partial \theta_i
\partial \theta_j }\big (
X_t-f_{\theta_0}^t \big)\Big]=\\
2\E\Big[\mbox{Tr}\Big
((H_{\theta_0}^t  \big )^{-2} \frac {\partial H_{\theta_0}^t  }{\partial \theta_j} \frac
{\partial H_{\theta_0}^t  }{\partial \theta_i}\Big )-\mbox{Tr}\Big
(\big(H_{\theta_0}^t  \big )^{-1}\frac {\partial^2 H_{\theta_0}^t  }
{\partial \theta_j\partial \theta_i}\Big )\Big]
\end{multline*}
>From equation (\ref{eq::der2}), we then derive the explicit expression
\begin{equation}\label{F0}
\big (F(\theta_0)\big )_{ij}=\E \Big [ 2\Big (\frac {\partial f^t_{\theta_0}}{\partial
\theta_j }\Big )' \big (H_{\theta_0}^t \big )^{-1}\frac
{\partial
f^t_{\theta_0}}{\partial \theta_i }+ \mbox{Tr} \Big
((H_{\theta_0}^t  \big )^{-2} \frac {\partial H_{\theta_0}^t  }{\partial \theta_j} \frac
{\partial H_{\theta_0}^t  }{\partial \theta_i}\Big ) \Big ].
\end{equation}
Under Assumption {\bf (Var)}, $F(\theta_0)$ is a positive definite $d\times d$
matrix. Indeed, for all $Y=(y_1,\ldots,y_d)\in\R^d$,
\begin{multline*}
Y'F(\theta_0)Y=\E\Big [ 2\Big (\sum_{1\le i\le d}y_i\frac {\partial f^t_{\theta_0}}{\partial
\theta_i }\Big )' \big (H_{\theta_0}^t \big )^{-1}\Big(\sum_{1\le i\le d}y_i\frac
{\partial
f^t_{\theta_0}}{\partial \theta_i }\Big)+\\
\mbox{Tr} \Big
((H_{\theta_0}^t  \big )^{-2} \Big(\sum_{1\le i\le d}y_i\frac {\partial H_{\theta_0}^t  }
{\partial \theta_i}\Big)^2\Big ) \Big ].
\end{multline*}
These two terms are nonnegative and at least one of them is positive under Assumption {\bf (Var)}.
Then $F(\theta_0)$ is an invertible matrix and there
exists $n$ large enough such that $F_n$ is an invertible
matrix. Moreover, (\ref{eq::tayl}) implies,
$$
n(\widehat \theta_n-\theta_0)=-2F_n^{-1}\Big(\frac {\partial
L_n(\widehat \theta_n)} {\partial \theta}-\frac {\partial L_n(
\theta_0)}{\partial \theta}\Big ).
$$
Therefore, if $\displaystyle \frac 1 { \sqrt{n}}\Big  \| \frac
{\partial L_n(\widehat \theta_n)}{\partial \theta} \Big
\|\limiteproban 0$, using Lemma \ref{lem::mart} one obtains
Theorem \ref{th::asnorm}. Since $\displaystyle \frac {\partial
\widehat L_n(\widehat \theta_n)}{\partial \theta} =0$ ($\widehat
\theta_n$ is a local extremum for $\widehat L_n$),
\begin{eqnarray}\label{L1}
\E \Big [ \frac 1 {\sqrt n} \Big \|  \frac {\partial L_n}{\partial
\theta} -\frac {\partial \widehat L_n}{\partial \theta} \Big
\|_\Theta \Big ] \limiten 0.
\end{eqnarray}
Using the relation (\ref{Dq}), the following inequality
$$|a_1b_1c_1-a_2b_2c_2| \leq
|a_1-a_2||b_2||c_2|+|a_1||b_1-b_2||c_2|+|a_1||b_1||c_1-c_2|$$
and the bounds $\| (\widehat  H_\theta^t )^{-1}\|_\Theta
\leq \underline H^{-1/m}$,
$ \| ( H_\theta^t  )^{-1}
\|_\Theta \leq \underline H^{-1/m}$, one obtains:
\begin{eqnarray*}
\Big \|  \frac {\partial q_t(\theta)}{\partial \theta_i} - \frac
{\partial \widehat q_t(\theta)}{\partial \theta_i}\Big \|_\Theta \hspace{-3mm}&
\leq &\hspace{-3mm} \frac 2 {\underline H^{1/m}}\Big[\Big  \|
\frac {\partial \widehat f^t_\theta}{\partial
\theta_i} -\frac {\partial f^t_\theta}{\partial \theta_i} \Big
\|_\Theta \big \| X_t-\widehat
f_\theta^t \big \|_\Theta \hspace{-1.5mm}+ \Big \|\frac {\partial f^t_\theta}{\partial
\theta_i} \Big \|_\Theta\big \| \widehat f_\theta^t-f_\theta^t \big \|_\Theta\Big]\\
&& \hspace{-4cm}
+2\Big \|\frac {\partial f^t_\theta}{\partial \theta_i} \Big
\|_\Theta \Big \| \big (H_\theta^t \big )^{-1}\hspace{-3mm}-\hspace{-1.5mm}\big (\widehat
H_\theta^t \big )^{-1}\Big \|_\Theta \big \| X_t-\widehat
f_\theta^t \big \|_\Theta \hspace{-1.5mm}+\big \| \widehat f_\theta^t-f_\theta^t \big
\|_\Theta \Big \|\frac {\partial \big (\widehat H_\theta^t \big
)^{-1}}{\partial \theta_i} \Big \|_\Theta \big \| X_t-\widehat
f_\theta^t \big \|_\Theta \\
&&  \hspace{-4cm}+\big \| X-f_\theta^t \big \|_\Theta \big
\| X_t-\widehat f_\theta^t \big \|_\Theta \Big \|\frac {\partial
\big (H_\theta^t \big )^{-1}}{\partial \theta_i}-\frac {\partial
\big (\widehat H_\theta^t  \big
)^{-1}}{\partial \theta_i} \Big \|_\Theta \hspace{-1.5mm} + \Big \| \big (\widehat H_\theta^t \big
)^{-1}\Big \|_\Theta \Big \|\frac {\partial H_\theta^t}{\partial \theta_i} -
\frac {\partial \widehat H_\theta^t
}{\partial \theta_i} \Big \|_\Theta\\
&& \hspace{2cm} +\Big \|\big (H_\theta^t \big )^{-1}- \big
(\widehat  H_\theta^t  \big )\Big \|_\Theta \Big \|\frac {\partial
\big (H_\theta^t \big )^{-1}}{\partial \theta_i} \Big \|_\Theta
\end{eqnarray*}
Under  {\bf (A3(f))} and {\bf (A3(M))} or {\bf (A3(H))}, there exists $C>0$ such that
$$\E \big \|
f^t_\theta - \widehat f^t_\theta \big \|_\Theta^r
\hspace{-1.5mm}\leq C \Big (\sum_{j \geq t} \alpha_j(f) \Big
)^r\mbox{and } \E \Big \|\frac {\partial
f^t_\theta}{\partial \theta_i} - \frac {\partial \widehat
f^t_\theta}{\partial \theta_i}  \Big \|_\Theta^r\hspace{-1.5mm}\leq
C \Big (\sum_{j \geq t} \alpha^{(1)}_j(f)
\Big )^r.
$$
The differences $\displaystyle\E \big \| H_{\theta}^t -\widehat
H_{\theta}^t  \big \|^{r/2}_\Theta  \leq C \Big (\sum_{j \geq t}
\alpha_j(M) \Big )^{r/2}$ can also be bounded:
\begin{eqnarray*}
\E \Big \|\frac {\partial
H_\theta^t }{\partial \theta_i} -\frac {\partial \widehat
H_\theta^t }{\partial \theta_i} \Big \|^{r/2}_\Theta&\hspace{-3mm}
\leq&\hspace{-3mm} C\Big (\Big (\sum_{j
\geq t} \alpha_j(M) \Big )^{r/2}\hspace{-3mm}+\hspace{-0.1cm} \Big (\sum_{j \geq t}
\alpha^{(1)}_j(M) \Big )^{r/2}\Big
),\\
\displaystyle\E \Big \|\frac {\partial \big (
H_\theta^t \big)^{-1}}{\partial \theta_i}
- \frac {\partial  \big ( \widehat H_\theta^t
\big)^{-1}}{\partial \theta_i} \Big \|^{r/2}_\Theta
\hspace{-0.1cm}&\hspace{-3mm}\leq & \hspace{-3mm}C  \Big (\Big (\sum_{j \geq t} \alpha_j(M) \Big
)^{r/2}\hspace{-0.4cm} + \hspace{-0.1cm} \Big (\sum_{j \geq t}
\alpha^{(1)}_j(M) \Big )^{r/2}\Big ).
\end{eqnarray*}
Finally, using H\"older inequalities, it exists another constant $C \geq 0$
satisfying
\begin{multline*}
\E\Big \|  \frac {\partial q_t(\theta)}{\partial \theta_i} - \frac
{\partial \widehat q_t(\theta)}{\partial \theta_i}\Big \|_\Theta
\leq  C \sum_{j \geq t} \big ( \alpha_j(f)+
\alpha_j(M)+\alpha_j(H)\\
+\alpha^{(1)}_j(f)+\alpha^{(1)}_j(M)+\alpha^{(1)}_j(H) \big ).
\end{multline*}
Under (\ref{condth2}), $\displaystyle \frac 1 {\sqrt n} \sum_{t=1}^n
 \E  \Big \| \frac {\partial q_t(\theta)}{\partial \theta_i} - \frac
{\partial \widehat q_t(\theta)}{\partial \theta_i}\Big \|_\Theta
\limiten  0$, and Theorem \ref{th::asnorm} follows.
\subsubsection*{Proof of Lemma \ref{lem::prel}}
Here, we focus on the case of $H_\theta$ under {\bf (A3(f))} and {\bf (A3(M))}.
The other cases are simpler.\\

With the same method and notation as in the proof of Lemma
\ref{lem::estvar}, the result holds as soon as the function
$\theta \in \Theta \to H^{t,p}_\theta$ is proved to satisfy a
Cauchy criterion in $\L^{r/2}\big ({\cal D}^{(2)}{\cal
C}(\Theta,{\cal M}_m)\big )$. Using the proof of Lemma
\ref{lem::estvar}, we already have
$\E\|H^{t,p}_\theta\|_\Theta^{r/2}<\infty$. It remains to bound
the quantities
$$
\E\Big
\|\frac {\partial H_\theta^{t,p}}{\partial \theta_i} \Big \|^{r/2}_\Theta\mbox{ and }
\E\Big
\|\frac {\partial^2 H_\theta^{t,p}}{\partial \theta_i\partial \theta_j}
\Big \|^{r/2}_\Theta\qquad\forall\,i, j\in \{1,\ldots,d\},\qquad\forall p\in\N^\ast.
$$
Using Assumption {\bf (A3(M))}:
\begin{multline*}
\Big\|\frac {\partial H_\theta^{t,p}}{\partial \theta_i} \Big
\|_\Theta\le 2\|M_\theta^{t,p}\|_\Theta\Big\|\frac {\partial
M_\theta^{t,p}}{\partial \theta_i} \Big \|_\Theta\\
\le\Big(\|M_\theta(0)\|_\Theta+\sum_{j=1}^{\infty}
\alpha_j(M)\|X_{t-j}\|\Big)\Big(\Big\|\frac {\partial M_\theta(0)}
{\partial \theta_i}\Big\|_\Theta+\sum_{j=1}^{\infty}
\alpha_j^{(1)}(M)\|X_{t-j}\|\Big).
\end{multline*}
Using $\E \big [ \|X_0\|^r\big ]<\infty$ and the H\"older and
Minkowsky inequalities:
\begin{multline*}
\E\Big [\|\frac {\partial H_\theta^{t,p}}{\partial \theta_i} \Big
\|^{r/2}_\Theta \Big ] \leq
C\Big(\|M_\theta(0)\|_\Theta^r+\E\big[\|X_0\|^r\big]\big(\sum_{j=1}^{\infty}
\alpha_j(M)\big)^r\Big)^{1/2}\\
\times\Big(\Big\|\frac {\partial M_\theta(0)}{\partial
\theta_i}\Big\|_\Theta^r+
\E\big[\|X_0\|^r\big]\big(\sum_{j=1}^{\infty}
\alpha_j^{(1)}(M)\big)^r\Big)^{1/2}.
\end{multline*}
In the same way, there exists another constant $C>0$ such that
\begin{align*}
\E\Big
\|\frac {\partial^2 H_\theta^{t,p}}{\partial \theta_i\partial \theta_j} \Big \|^{r/2}_\Theta\le&
C\Big[\Big(\big(\sum_{j=1}^{\infty}
\alpha_j^{(1)}(M)\big)^r\big(\sum_{j=1}^{\infty}
\alpha_j^{(1)}(M)\big)^r\Big)^{1/2}\\
&+\Big(\big(\sum_{j=1}^{\infty}
\alpha_j(M)\big)^r\big(\sum_{j=1}^{\infty}
\alpha_j^{(2)}(M)\big)^r\Big)^{1/2}\Big].
\end{align*}
>From $\sum_j\alpha_j(M)\hspace{-1mm}<\hspace{-1mm}\infty $,
$\sum_j\alpha^{(1)}_j(M)\hspace{-1mm}<\hspace{-1mm}\infty$ and
$\sum_j\alpha^{(2)}_j(M)\hspace{-1mm}<\hspace{-1mm}\infty$ we
deduce that $\E \big [\|H^{t,p}_\theta\|_{(2),\Theta}^{r/2} \big
]\hspace{-1mm} <\hspace{-1mm} \infty$ for all $p\in\N^\ast$. In
the same way as in the proof of Lemma \ref{lem::estvar} we can
also prove that the sequence $(H_\theta^{t,p})_{p\in\N^\ast}$
satisfies the Cauchy criterion in the Banach space $\L^{r/2}({\cal
D}^{(2)}{\cal C}(\Theta,{\cal M}_m))$. For the first derivatives,
the result easily follows from the inequality
$$
\Big\|\frac {\partial H_\theta^{t,p}}{\partial \theta_i} - \frac
{\partial H_\theta^{t,q}}{\partial \theta_i}\Big
\|_\Theta\hspace{-1mm}\le
2\|M_\theta^{t,p}-M_\theta^{t,q}\|_\Theta\Big\|\frac {\partial
M_\theta^{t,p}} {\partial \theta_i} \Big
\|_\Theta\hspace{-1mm}+2\|M_\theta^{t,q}\|_\Theta \Big\|\frac
{\partial H_\theta^{t,p}}{\partial \theta_i} -\frac {\partial
H_\theta^{t,q}}{\partial \theta_i}\Big \|_\Theta\hspace{-1mm}.
$$
For the second derivatives, a similar argument finishes the proof.
\subsubsection*{Proof of Lemma \ref{lem::mart}}
Simple calculations give the relations
$$
\frac {\partial \big (H_\theta^t  \big )^{-1}}{\partial \theta_k}=
-\big (H_\theta^t  \big )^{-1}\frac {\partial H_\theta^t}
{\partial \theta_k}\big (H_\theta^t  \big )^{-1}\mbox{ and }
\frac {\partial \ln\det \big( H_\theta^t  \big) }
{\partial \theta_k}=\mbox{Tr} \Big
( \big (H_\theta^t  \big )^{-1} \frac {\partial H_\theta^t}{\partial \theta_k} \Big ).
$$
>From Lemma \ref{lem::prel}, ${\partial f_\theta^t}/{\partial
\theta}$, ${\partial H_\theta^t}/{\partial \theta}$ and $\big
(\widehat H_\theta^t  \big) ^{-1}$ are a.s. finite. Then
${\partial L_n(\theta)}/{\partial \theta}$ is an a.s. finite
measurable function satisfying, for all $1\leq i \leq d$,
${\partial L_n(\theta)}/{\partial \theta_i}=-\frac 1 2
\sum_{t=1}^n{\partial q_t(\theta)}/{\partial \theta_i}$ with
\begin{multline}
\label{Dq} \frac {\partial q_t(\theta)}{\partial \theta_k} =-2
\Big (\frac {\partial f^t_\theta}{\partial \theta_k}\Big )' \big
(H_\theta^t  \big )^{-1} \big ( X_t-f_\theta^t \big) \\
+ \big (
X_t-f_\theta^t \big)' \frac {\partial \big (H_\theta^t \big
)^{-1}}{\partial \theta_k} \big (
X_t-f_\theta^t \big)+ \mbox{Tr} \Big
( \big (H_\theta^t  \big )^{-1} \frac {\partial H_\theta^t}{\partial \theta_k} \Big ).
\end{multline}
Denoting ${\mathcal F}_t=\sigma(X_t,X_{t-1},\ldots)$, let us prove
that $\displaystyle \Big ( \frac {\partial q_t(\theta_0)}{\partial
\theta},{\mathcal F}_t\Big )_{t\in\Z}$ is a $\R^m$-valued
martingale difference process. Indeed, for all $t \in \Z$,
$$
\E \big ( (X_t-f_{\theta_0}^t)
|{\mathcal
F}_t)=0\mbox{ and }\E\big( \big (
X_t-f_{\theta_0}^t \big) \big ( X_t-f_{\theta_0}^t \big)'|{\mathcal
F}_t\big)=H_{\theta_0}^t.
$$
As a consequence,
$$
\E \Big ( \frac {\partial q_t(\theta_0)}{\partial
\theta_k}|{\mathcal
F}_t\Big )\hspace{-1mm}=\hspace{-1mm} \E \Big ( \big ( X_t-f_{\theta_0}^t \big)'
\frac {\partial \big (H_{\theta_0}^t  \big )^{-1}}{\partial \theta_k} \big (
X_t-f_{\theta_0}^t \big)|{\mathcal
F}_t\Big )+\mbox{Tr} \Big ( \big
(H_{\theta_0}^t  \big )^{-1} \frac {\partial H_{\theta_0}^t}{\partial \theta_k} \Big
).
$$
We conclude by noticing that the first term of the sum is equal to
$$\E \Big (\mbox{Tr} \Big (\frac
{\partial  \big (H_{\theta_0}^t \big )^{-1} }{\partial \theta_k}
 \big ( X_t-f_{\theta_0}^t \big)\big (
X_t-f_{\theta_0}^t \big)'\Big )|{\mathcal F}_t\Big )=\mbox{Tr}
\Big (\frac {\partial\big (H_{\theta_0}^t \big )^{-1} }{\partial
\theta_k}H_{\theta_0}^t\Big ).$$ In order to apply the Central
Limit Theorem for martingale-differences, see \cite{Billingsley1968},
we have to prove that $\displaystyle \E \Big [ \Big \|\frac
{\partial q_t(\theta_0)}{\partial \theta} \Big \|^2\Big ]
<\infty$. Using the relation
$X_t-f_{\theta_0}^t=M_{\theta_0}^t\xi_t$ for all $t\in\Z$, then
\begin{align*}
\frac {\partial q_t(\theta_0)}{\partial \theta_k}= &-2\Big
(\frac {\partial f^t_{\theta_0}}{\partial \theta_k}\Big )'\big
(H_{\theta_0}^t \big )^{-1} M_{\theta_0}^t\xi_t-\xi_t'{M_{\theta_0}^t}'{\big
(H_{\theta_0}^t \big )^{-1}}'  \frac
{\partial H_{\theta_0}^t}{\partial \theta_k}
\big
(H_{\theta_0}^t \big )^{-1} M_{\theta_0}^t\xi_t\\
&+\mbox{Tr} \Big (
\big (H_{\theta_0}^t \big )^{-1} \frac {\partial H_{\theta_0}^t }{\partial
\theta_k} \Big )
\end{align*}
Let us compute the expectation of the square of the second term of
the sum, with $\mbox{Tr} (ABC)=\mbox{Tr}(CAB)=\mbox{Tr}(ACB)$ for
symmetric matrices $A$, $B$ and $C$,
\begin{align}\label{calcder}
&\E\Big[(\xi_t'\xi_t)^2\mbox{ Tr}\Big({M_{\theta_0}^t}'{\big
(H_{\theta_0}^t \big )^{-1}}  \frac
{\partial H_{\theta_0}^t}{\partial \theta_k}
{\big
(H_{\theta_0}^t \big )^{-1}} \frac
{\partial H_{\theta_0}^t}{\partial \theta_k}\big
(H_{\theta_0}^t \big )^{-1}M_{\theta_0}^t\Big)\Big]\\
\nonumber &=\E\Big[(\xi_t'\xi_t)^2\mbox{ Tr}\Big({\big
(H_{\theta_0}^t \big )^{-2}}  \Big(\frac
{\partial H_{\theta_0}^t}{\partial \theta_k}
 \Big)^2\Big)\Big].
\end{align}
Using this relation, the bound $  \|\big(H_{\theta_0}^t \big)^{-1}\|_\Theta \leq \underline H
^{-1/m}$ and the independence of $\xi_t$ and ${\cal F}_t$, there exists $C>0$ such that
\begin{align*}\E \Big [\Big (\frac {\partial
q_t(\theta_0)}{\partial \theta_k} \Big )^2 \Big ] \le& C \left (
\E \Big [\Big \| \frac {\partial f^t_{\theta_0}}{\partial
\theta_k} \Big \|^2\Big \| M_{\theta_0}^t \Big \|^2\Big ] \times \E \big
[ \| \xi_t \|^2\big ]+\E \Big [ \Big \| \frac {\partial H_{\theta_0}^t }
{\partial \theta_k} \Big \|^2 \Big ] \right .\\
&\left . + \E \big [ \| \xi_t' \xi_t \|^2 \big ] \times
\E \Big [ \Big \| \frac {\partial H_{\theta_0}^t  }{\partial \theta_k} \Big \|^2 \Big ] \right )
\end{align*}
Therefore, since $r\ge 4$, the moment conditions for the CLT are fulfilled
$$
\E \Big [\Big \| \frac {\partial q_t(\theta_0)}{\partial \theta}
\Big \|^2 \Big ]=\sum_{k=1}^d \E \Big [\frac {\partial
q_t(\theta_0)}{\partial \theta_k} \Big ]^2< \infty.
$$
We compute the asymptotic covariance matrix of $\displaystyle \frac
{\partial q_t(\theta_0)}{\partial \theta}$.
Thus, $(G(\theta_0))_{ij}$
\begin{eqnarray}
\nonumber &=&
\E \Big [ \frac {\partial
q_t(\theta_0)}{\partial \theta_i}\frac {\partial
q_t(\theta_0)}{\partial \theta_j} \Big ]\\
\nonumber &=&
\E \left [4\Big (\frac {\partial
f^t_{\theta_0}}{\partial \theta_i}\Big )'\big (H_{\theta_0}^t
\big )^{-1}\Big (\frac {\partial f^t_{\theta_0}}{\partial
\theta_j}\Big )-\mbox{Tr}\Big ( \big(H_{\theta_0}^t \big )^{-1}\frac
{\partial H_{\theta_0}^t}{\partial \theta_i}\Big)\mbox{Tr}\Big (  \big
(H_{\theta_0}^t \big )^{-1}\frac {\partial H_{\theta_0}^t}{\partial \theta_j}\Big)
 \right . \\
\label{G0}  &&
\left . + p\big (m_4+(p-1)\big
)\mbox{Tr} \Big (\big (H_{\theta_0}^t \big
)^{-2} \frac {\partial H_{\theta_0}^t}{\partial \theta_i}
\frac {\partial H_{\theta_0}^t}{\partial \theta_j}\Big )  \right ].
\end{eqnarray}
To simplify the expression, we assume here that $\xi_t$ and $-\xi_t$
have the same distribution in order that $\E \big [ \xi_t \xi_t' A
\xi_t \big ]=0$ for $A$ a matrix.
\subsubsection*{Proof of the Lemma \ref{lem::der2}}
>From the proof of Proposition \ref{stationarity} and from the result
of Lemma \ref{lem::prel}, the second derivative process
$({\partial^2 q_t(\theta)}/{\partial \theta^2})_{t\in\Z}$ is
stationary ergodic (it is a measurable function of
$X_t,X_{t-1},\ldots$). Therefore it satisfies a Uniform Law of Large Numbers (ULLN) if its first
uniform moment is bounded. \\

>From equation (\ref{Dq}), the second partial derivatives of
$q_t(\theta)$ are
\begin{eqnarray}
\nonumber\frac {\partial^2 q_t(\theta)}{\partial \theta_i \partial
\theta_j} \hspace{-3mm}&=&\hspace{-3mm} -2 \Big (\frac {\partial^2 f^t_\theta}{\partial
\theta_i \partial \theta_j}\Big )' \big (H_\theta^t \big )^{-1}
\big ( X_t-f_\theta^t \big)\\
\nonumber&& \hspace{-6mm}+ \big ( X_t-f_\theta^t \big)' \frac
{\partial^2 \big (H_\theta^t  \big )^{-1}}{\partial \theta_i
\partial \theta_j} \big ( X_t-f_\theta^t \big) -2 \Big (\Big (\frac {\partial
f^t_\theta}{\partial \theta_i }\Big )' \frac {\partial \big
(H_\theta^t  \big )^{-1}}{\partial \theta_j}\\
\nonumber&& \hspace{-6mm} + \Big (\frac
{\partial f^t_\theta}{\partial \theta_j }\Big )' \frac {\partial
\big (H_\theta^t  \big )^{-1}}{\partial \theta_i} \Big )\big (
X_t-f_\theta^t \big)+ 2\Big (\frac {\partial f^t_\theta}{\partial
\theta_i }\Big )' \big (H_\theta^t \big )^{-1}\Big (\frac
{\partial f^t_\theta}{\partial \theta_i }\Big )   \\
\label{eq::der2}&&\hspace{-6mm} +\mbox{Tr} \Big (\Big ( \frac {\partial \big
(H_\theta^t \big )^{-1} }{\partial \theta_j}\Big )\Big ( \frac
{\partial H_\theta^t  }{\partial \theta_i}\Big )\Big
) +\mbox{Tr} \Big (\big (H_\theta^t  \big )^{-1} \Big ( \frac
{\partial^2 H_\theta^t  }{\partial \theta_i \partial
\theta_j}\Big )\Big ).
\end{eqnarray}
Therefore, using the
bound $\big \|\big (H_\theta^t  \big )^{-1}\big \|_\Theta \leq
\underline{M}^{-1/m}$ of Lemma \ref{lem::estvar} and usual
relations between norms and traces of matrix, there exists $C>0$ such that
\begin{multline*}
\Big \| \frac {\partial^2 q_t({\theta})}{\partial \theta_i
\partial \theta_j}\Big \|_\Theta \leq C \Big [\Big (
\Big \| \frac {\partial^2 f^t_{\theta}}{\partial \theta_i
\partial \theta_j}\Big \|_\Theta + \Big \| \frac {\partial H_{\theta}^t }
{\partial \theta_j} \Big \|_\Theta  \Big \| \frac {\partial
f^t_{\theta}}{\partial \theta_i }\Big \|_\Theta +\Big \| \frac
{\partial H_{\theta}^t }{\partial \theta_i} \Big \|_\Theta
\Big \| \frac {\partial
f^t_{\theta}}{\partial \theta_j }\Big \| \Big ) \big \|X_t-f^t_{\theta}\big \|_\Theta\\
+ \Big \| \frac {\partial^2 H_{\theta}^t }{\partial
\theta_i
\partial \theta_j} \Big \|_\Theta\big \| X_t-f^t_{\theta}
 \big \|_\Theta^2
+ \Big \| \frac {\partial f^t_{\theta}}{\partial \theta_i }\Big
\|_\Theta \Big \| \frac {\partial f^t_{\theta}}{\partial \theta_j }\Big
\|_\Theta + \Big \| \frac {\partial H_\theta^t  }{\partial \theta_i}
\Big \|_\Theta \Big \| \frac {\partial H_\theta^t  }{\partial \theta_j} \Big \|_\Theta\Big ].
\end{multline*}
We conclude that $\displaystyle\E\Big \| \frac {\partial^2 q_t({\theta})}{\partial \theta_i
\partial \theta_j}\Big \|_\Theta^{r/4}<\infty$ ($r\ge 4$) since, for $t \in \Z$, $1\leq i,j \leq d$,
\begin{gather*}
\E \big [ \| X_t\|^r \big ]<+\infty, \E \big [\big \| f_\theta^t \big \|^r_\Theta \big ]
<+\infty, \E \Big [ \Big \|\frac {\partial
f^t_\theta}{\partial \theta_i} \Big \|_\Theta^r \Big ] <+\infty,
\E \Big [ \Big \|\frac {\partial^2
f^t_\theta}{\partial \theta_i
\partial \theta_j}\Big \|_\Theta^r \Big ]
<+\infty;\\
 \E
\Big [ \Big \| H_{\theta}^t  \Big \|_\Theta ^{r/2} \Big ]<+\infty,
\E \Big [ \Big \| \frac {\partial H_{\theta}^t}{\partial \theta_i}
 \Big \|_\Theta ^{r/2}
\Big ]<+\infty, \displaystyle \E \Big [ \Big \| \frac {\partial^2
H_{\theta}^t}{\partial \theta_i
\partial \theta_j} \Big \|_\Theta ^{r/2} \Big ]<\infty.
\end{gather*}
As a consequence, the ULLN holds for ${\partial^2
q_t({\theta})}/{\partial \theta^2}$.\\

{\bf Acknowledgements.} We would are very grateful to Thomas
Mikosch who made a critical review of the drafts and with whom we
have worked on the final version of this paper.
\bibliographystyle{acm}%plain}%
\bibliography{BardetWintenberger}

\end{document}